\documentclass{article}

\usepackage{ifthen}
\newboolean{neurips_version}
\setboolean{neurips_version}{true}
\ifthenelse{\boolean{neurips_version}}
{\PassOptionsToPackage{numbers, compress}{natbib}
 \usepackage[preprint]{neurips_2022}
}
{
\usepackage{PRIMEarxiv}
 \usepackage[numbers]{natbib}
}

\usepackage[utf8]{inputenc} 
\usepackage[T1]{fontenc}    
\usepackage{hyperref}       
\usepackage{url}            
\usepackage{booktabs}       
\usepackage{amsfonts,amsmath,amsthm,amssymb}       
\usepackage{nicefrac}       
\usepackage{microtype}      
\usepackage{xcolor}         
\usepackage[disable]{todonotes}
\usepackage{enumitem}
\usepackage{caption}

\usepackage{mathtools} 
\usepackage{tikz}
\usetikzlibrary{calc}
\usepackage{pgfplots}
\usepackage{graphicx}

\usepackage{amsmath}
\usepackage{algorithm, algpseudocode}
\usepackage{amsfonts}  

\usepackage[english]{babel}
\usepackage{amsthm}

\usepackage{subfig}
\usepackage[export]{adjustbox}

\usepackage{hhline}
\usepackage{makecell, caption, booktabs}
\usepackage{siunitx}

\usepackage{multirow}

\usepackage{amsmath, nccmath}
\usepackage{bigstrut}

\usepackage{booktabs}

\usepackage{lipsum}
\graphicspath{{media/}}     

\usepackage{CJKutf8}
\usepackage[framed,numbered,autolinebreaks,useliterate]{mcode}

\title{ Generalized Quadratic Gradient: A New Direction in Optimization via the Fusion of Positive-Definite Curvature Matrices and Gradients into A Unified Framework }

\author{ \href{https://orcid.org/0000-0003-0378-0607}{\includegraphics[scale=0.06]{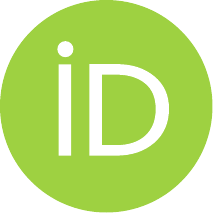}\hspace{1mm}John Chiang} \\                             
                                      \\
	\texttt{john.chiang.smith@gmail.com} 
}

\date{}

\newtheorem*{definition*}{Definition} 

\theoremstyle{remark}

\renewcommand{\epsilon}{\varepsilon}

\makeatletter
\def\namedlabel#1#2{\begingroup
    #2%
    \def\@currentlabel{#2}%
    \phantomsection\label{#1}\endgroup
}
\makeatother

\algnewcommand{\LeftComment}[1]{\Statex \(\triangleright\) #1}
\algnewcommand{\LineCommentStep}[1]{\Statex \textbf{[Step #1]:} }
\makeatletter
\newlength{\trianglerightwidth}
\algnewcommand{\LineComment}[1]{\Statex \hskip\ALG@thistlm $\triangleright$ #1}
\algnewcommand{\LineCommentCont}[1]{\Statex \hskip\ALG@thistlm%
  \parbox[t]{\dimexpr\linewidth-\ALG@thistlm}{\hangindent=\trianglerightwidth \hangafter=1 \strut$\triangleright$ #1\strut}}
\algnewcommand{\LeftLineCommentCont}[1]{\Statex \hskip\ALG@thistlm%
  \parbox[t]{\dimexpr\linewidth-\ALG@thistlm}{\leftskip=\algorithmicindent \hangindent=\trianglerightwidth \hangafter=1 \strut$\triangleright$ #1\strut}}

\begin{document}

\maketitle

\begin{abstract}%

Quadratic Gradient (QG) is a Newton-type optimization framework that bridges first-order gradient descent and second-order optimization by incorporating curvature information into gradient updates. Simplified Quadratic Gradient (SQG) reduces the complexity of QG construction while preserving its optimization capability, whereas Quasi-Quadratic Gradient (QQG) extends the quadratic gradient principle to quasi-Newton methods such as BFGS.

In this paper, we propose **Generalized Quadratic Gradient (GQG)**, a unified framework that extends the quadratic gradient principle to a broader class of Newton-type optimization algorithms. By abstracting the common structure of existing quadratic gradient methods, we show that the fundamental requirement of quadratic gradient construction is not limited to specific Hessian approximations, such as constant Hessian matrices, diagonal Hessian approximations, or BFGS-based Hessian surrogates. Instead, it can be generalized to any positive-definite curvature matrix satisfying the stationary condition of a local quadratic model.

Based on this perspective, we investigate the construction of generalized quadratic gradients using various positive-definite Hessian surrogates beyond BFGS, providing a broader foundation for developing curvature-aware optimization algorithms.

\end{abstract}

\listoftodos

\section{Introduction}

\subsection{Background}

Quadratic Gradient (QG) is a Newton-type optimization framework that bridges first-order gradient descent methods and second-order optimization algorithms. By incorporating curvature information into gradient updates, QG constructs a quadratic gradient direction that captures the benefits of Newton optimization while avoiding the explicit computation of the exact Hessian matrix.

Simplified Quadratic Gradient (SQG) further reduces the complexity of the original quadratic gradient construction by simplifying the Hessian approximation, while maintaining comparable optimization performance. Subsequently, Quasi-Quadratic Gradient (QQG) extends the quadratic gradient principle beyond fixed Hessian-based methods by incorporating quasi-Newton techniques, such as BFGS, into the quadratic gradient framework.

In this paper, we propose **Generalized Quadratic Gradient (GQG)**, a general framework that investigates how the quadratic gradient principle can be extended to a broader class of Newton-type second-order optimization algorithms. By abstracting the common principles behind QG, SQG, and QQG, GQG provides a generalized formulation for constructing quadratic-gradient-based optimization methods.

Unlike previous quadratic gradient formulations that mainly rely on constant (diagonal) Hessian approximations derived from (Simplified) Fixed Hessian methods or Hessian surrogates generated by BFGS, we show that the fundamental requirement of quadratic gradient construction should not be restricted to a specific Hessian approximation strategy. Instead, the essential requirement is the construction of a positive-definite curvature matrix that satisfies the stationary condition of a local quadratic model.

Therefore, we investigate various positive-definite Hessian surrogate construction strategies beyond BFGS for developing generalized quadratic gradient algorithms. Specifically, given a positive-definite curvature matrix, GQG defines the optimization direction by minimizing a local quadratic model:

\[
m(p)=g^Tp+\frac{1}{2}p^TBp,
\]

where ($B\succ0$) represents a generalized curvature matrix. The stationary condition of this quadratic model yields:

\[
p=-B^{-1}g,
\]

which provides a unified formulation connecting gradient descent and Newton-type optimization methods.

Finally, we evaluate the performance of different generalized quadratic gradient variants constructed from various positive-definite curvature approximations. Extensive experiments are conducted to investigate whether the quadratic gradient principle can be effectively generalized beyond fixed Hessian and BFGS-based approaches.


Gradient-based optimization methods are fundamental tools for solving large-scale machine learning and numerical optimization problems. Among them, first-order methods such as gradient descent achieve remarkable scalability by relying only on gradient information, but they often suffer from slow convergence due to the lack of curvature awareness. In contrast, second-order methods, particularly Newton-type algorithms, exploit the local curvature of the objective function through the Hessian matrix and can achieve significantly faster convergence near the optimum. However, the computational and memory costs associated with Hessian construction and inversion limit their applicability to large-scale problems.

A common strategy for bridging the gap between first-order and second-order optimization is to approximate the curvature information instead of computing the exact Hessian. This leads to a broad family of curvature-aware optimization methods, including quasi-Newton methods, Gauss-Newton methods, natural gradient methods, and adaptive second-order optimization algorithms. Although these methods employ different Hessian approximation strategies, they share a common principle: constructing an appropriate curvature matrix to define a better local geometry for gradient updates.

Quadratic Gradient (QG) provides another perspective for combining gradient-based and Newton-type optimization. Instead of directly updating parameters using the gradient direction, QG derives the update direction by minimizing a local quadratic model:

\[
m(p)=g^Tp+\frac{1}{2}p^TBp,
\]

where \(g\) denotes the gradient and \(B\) represents a curvature matrix. The stationary condition of this quadratic model yields:

\[
p=-B^{-1}g,
\]

which recovers the Newton update when \(B\) is the exact Hessian and provides a generalized curvature-aware gradient update when \(B\) is an approximation.

The original Quadratic Gradient framework constructs the curvature matrix based on fixed Hessian information. Simplified Quadratic Gradient (SQG) further reduces computational complexity by simplifying the Hessian construction while maintaining similar optimization behavior. Subsequently, Quasi-Quadratic Gradient (QQG) extends the quadratic gradient principle to quasi-Newton optimization by incorporating BFGS-based Hessian surrogates.

However, existing quadratic gradient formulations remain closely tied to specific curvature construction strategies. This raises an important question: **is the quadratic gradient principle fundamentally dependent on a particular Hessian approximation method, or can it be generalized to a broader class of curvature matrices?**

In this work, we answer this question by introducing **Generalized Quadratic Gradient (GQG)**, a unified framework that extends the quadratic gradient principle beyond fixed Hessian and BFGS-based approaches. We show that the essential requirement for constructing a quadratic gradient is not a specific approximation of the Hessian matrix, but rather a positive-definite curvature matrix that defines a valid local quadratic geometry.

Based on this observation, GQG explores the construction of quadratic gradients using general positive-definite curvature matrices beyond BFGS. This provides a broader perspective for understanding quadratic-gradient-based optimization and establishes a flexible framework for developing new Newton-type optimization algorithms.

\subsection{Related Work}

\subsubsection{Gradient-Based Optimization}

First-order optimization methods are widely used in large-scale machine learning due to
their computational efficiency. The classical gradient descent (GD) updates model
parameters according to the negative gradient direction:
\[
x_{k+1}=x_k-\eta_k\nabla f(x_k),
\]
where $\eta_k$ denotes the learning rate. Although GD avoids the expensive computation
of second-order information, its convergence can be slow for ill-conditioned optimization
problems because it ignores the curvature of the objective landscape.

To improve convergence, momentum-based methods introduce historical gradient
information into the update process. Nesterov's accelerated gradient (NAG) achieves
accelerated convergence by evaluating gradients at an extrapolated point. Furthermore,
adaptive gradient methods, including AdaGrad and Adam, dynamically adjust the learning
rate using historical gradient statistics. These approaches have achieved remarkable
success in deep learning, but their curvature modeling is generally restricted to
diagonal gradient-based scaling rather than explicit second-order information.

\subsubsection{Newton-Type and Fixed-Hessian Optimization}

Newton-type methods improve upon first-order approaches by incorporating the local
curvature of the objective function through the Hessian matrix. Given a local quadratic
approximation:
\[
m(p)=f(x_k)+g_k^Tp+\frac{1}{2}p^TH_kp,
\]
the Newton direction is obtained by minimizing the quadratic model:
\[
p_k=-H_k^{-1}g_k.
\]

Although Newton methods provide fast local convergence, computing and inverting the
full Hessian matrix introduces significant computational and memory overhead.
Therefore, various Hessian approximation strategies have been investigated.

Fixed-Hessian Newton methods replace the iteration-dependent Hessian with a constant
positive-definite approximation. This strategy significantly reduces computational
cost while preserving the curvature-aware property of Newton optimization. Such
approaches have also been explored in privacy-preserving machine learning, where
reducing the number of optimization iterations is critical due to the high computational
cost of encrypted operations.

\subsubsection{Quasi-Newton Methods}

Quasi-Newton methods approximate second-order information without explicitly computing
the Hessian matrix. Instead, they construct a Hessian approximation using gradient
differences. Given
\[
s_k=x_{k+1}-x_k,\qquad y_k=g_{k+1}-g_k,
\]
the approximation matrix is updated according to the secant condition:
\[
B_{k+1}s_k=y_k.
\]

Among quasi-Newton algorithms, BFGS is one of the most successful approaches due to
its robustness and ability to maintain positive definiteness. Specifically, if the
initial approximation is positive definite and the curvature condition
$s_k^Ty_k>0$ is satisfied, the BFGS update preserves the positive definiteness of the
Hessian approximation, ensuring that the search direction remains a descent direction.

However, quasi-Newton methods rely on specific update rules to approximate the Hessian,
which limits their direct integration with other curvature-based optimization
frameworks.

\subsubsection{Adaptive Second-Order Optimization}

Recent studies have explored efficient approximations of second-order information for
large-scale optimization. Instead of computing the full Hessian matrix, these methods
construct simplified curvature representations.

AdaHessian estimates diagonal Hessian information to improve adaptive optimization
while maintaining computational efficiency. Other approaches, such as Shampoo and
K-FAC, exploit structured matrix approximations, including Kronecker-factorized
curvature matrices, to capture richer second-order information.

Although these methods employ different approximation strategies, they share a common
objective: constructing an efficient curvature matrix that improves the geometry of
gradient-based optimization.

\subsubsection{Quadratic Gradient Optimization}

Quadratic Gradient (QG) was introduced as a framework that bridges first-order gradient
methods and Newton-type optimization by incorporating curvature information into the
gradient update. Instead of directly following the gradient direction, QG constructs a
quadratic-gradient direction based on a Hessian proxy, which serves as an efficient
approximation of the inverse Hessian.

Simplified Quadratic Gradient (SQG) further reduces the complexity of QG construction
by replacing the original curvature approximation with a simplified positive-definite
proxy, demonstrating that effective second-order acceleration does not necessarily
require accurate Hessian estimation.

More recently, Quasi-Quadratic Gradient (QQG) extended the quadratic gradient principle
to quasi-Newton optimization by incorporating BFGS-based Hessian approximations.
This extension demonstrated that quadratic-gradient updates can be combined with
adaptive curvature estimation beyond fixed-Hessian methods.

However, existing quadratic-gradient formulations remain dependent on specific
curvature construction strategies, including fixed Hessian approximations, simplified
diagonal approximations, and BFGS-based Hessian surrogates. This motivates a more
general investigation of the fundamental requirement behind quadratic-gradient
optimization.

\subsection{Contributions}
The contributions of this work are summarized as follows:

1. We introduce Generalized Quadratic Gradient (GQG), a unified framework that extends the quadratic gradient principle from fixed Hessian and quasi-Newton methods to a broader family of Newton-type optimization algorithms.

2. We reveal that the essential requirement of quadratic gradient optimization is not a specific Hessian approximation technique, but rather a positive-definite curvature matrix satisfying the stationary condition of a local quadratic model.

3. We investigate various positive-definite curvature construction strategies beyond BFGS and develop corresponding generalized quadratic gradient algorithms.

4. We experimentally evaluate the effectiveness of these generalized quadratic gradient variants and analyze the impact of different curvature approximations on optimization performance.

\begin{itemize}

\item We identify the underlying mathematical principle of quadratic-gradient
optimization as positive-definite curvature transformation.

\item We propose the Generalized Quadratic Gradient (GQG) framework, which
unifies various curvature approximation strategies under a common
Newton-type gradient formulation.

\item We demonstrate that existing curvature constructions, including fixed
Hessian, diagonal approximation, and BFGS-based quasi-Newton methods, can be
naturally integrated into modern optimizers such as NAG, AdaGrad, and Adam.

\end{itemize}

In this work, we show that the quadratic gradient principle is not fundamentally
restricted to a particular Hessian approximation. Instead, it can be generalized to
any positive-definite curvature matrix satisfying the stationary condition of a local
quadratic model.

\section{Preliminaries}

\paragraph{Loewner Ordering}

For two symmetric matrices $A$ and $B$, the Loewner ordering is defined as:

\[
A \leq B
\]

if and only if

\[
B-A\succeq0,
\]

where $B-A$ is positive semi-definite. Equivalently, for any vector $\mathbf{x}$,

\[
\mathbf{x}^{T}A\mathbf{x}
\leq
\mathbf{x}^{T}B\mathbf{x}.
\]

This ordering is commonly used to define upper and lower bounds of Hessian matrices in
Newton-type optimization methods.

\subsection{Newton-Raphson method}

The Newton-Raphson method is a classical second-order optimization algorithm that
utilizes the Hessian matrix to capture the local curvature of the objective function.
Given a differentiable function $f(\mathbf{x})$, the parameter update is defined as:

\begin{equation}
    \mathbf{x}_{k+1}
    =
    \mathbf{x}_{k}
    -
    [\nabla^2 f(\mathbf{x}_k)]^{-1}
    \nabla f(\mathbf{x}_k),
\end{equation}

where $\nabla f(\mathbf{x}_k)$ and $\nabla^2 f(\mathbf{x}_k)$ denote the gradient and
Hessian matrix at iteration $k$, respectively. Although Newton-Raphson achieves
quadratic convergence near the optimum, explicitly computing and inverting the Hessian
matrix introduces substantial computational and memory costs, limiting its applicability
to high-dimensional optimization problems.

\subsubsection{Simplified Fixed Hessian}
To reduce the computational overhead of Newton-Raphson iterations, the Simplified Fixed
Hessian (SFH) method replaces the iteration-dependent Hessian matrix with a fixed
positive-definite approximation:

\begin{equation}
    \mathbf{x}_{k+1}
    =
    \mathbf{x}_{k}
    -
    \mathbf{H}_{fixed}^{-1}
    \nabla f(\mathbf{x}_k).
\end{equation}

By avoiding repeated Hessian evaluation and matrix inversion, SFH significantly reduces
the computational complexity while retaining the curvature-aware property of Newton-type
optimization. However, since the fixed Hessian approximation does not adapt to changes
in the optimization landscape, its convergence behavior depends strongly on the quality
of the selected matrix approximation, particularly for non-convex optimization problems.

\subsubsection{Quasi-Newton Optimization}

Quasi-Newton methods approximate second-order curvature information without explicitly
computing the Hessian matrix. Instead of evaluating
$\nabla^2 f(x_k)$ directly, they maintain a sequence of positive-definite matrices
that approximate either the Hessian or its inverse using only first-order gradient
information.

Compared with the classical Newton method, which requires Hessian evaluation and
matrix inversion with cubic computational complexity, quasi-Newton methods achieve a
more favorable computational cost while preserving fast convergence properties.

In the context of GQG, the inverse Hessian approximation generated by quasi-Newton
methods can be directly interpreted as the curvature transformation matrix:

\begin{equation}
    P_k \approx (\nabla^2 f(x_k))^{-1}.
\end{equation}

Therefore, quasi-Newton algorithms provide a natural family of curvature operators
for constructing generalized quadratic gradients.

\paragraph{BFGS Update Scheme}

Among quasi-Newton algorithms, the Broyden-Fletcher-Goldfarb-Shanno (BFGS) method is
one of the most widely adopted approaches due to its robustness and strong empirical
performance.

Given the step vector:

\begin{equation}
    s_k=x_{k+1}-x_k,
\end{equation}

and the gradient difference:

\begin{equation}
    y_k=g_{k+1}-g_k,
\end{equation}

the BFGS Hessian approximation update is defined as:

\begin{equation}
    B_{k+1}
    =
    B_k
    +
    \frac{y_ky_k^T}{y_k^Ts_k}
    -
    \frac{B_ks_ks_k^TB_k}{s_k^TB_ks_k}.
\end{equation}

The corresponding inverse Hessian approximation can be used as the curvature
transformation matrix in GQG:

\begin{equation}
    P_k=B_k^{-1}.
\end{equation}

A fundamental property of BFGS is that the positive definiteness of the Hessian
approximation is preserved. Specifically, if the initial matrix $B_0$ is symmetric
positive definite and the curvature condition

\begin{equation}
    s_k^Ty_k>0
\end{equation}

is satisfied, all subsequent matrices $B_k$ remain symmetric positive definite.

Consequently, the resulting inverse Hessian approximation also satisfies:

\begin{equation}
    B_k^{-1}\succ0,
\end{equation}

which fulfills the requirement of the GQG framework. This observation establishes
QQG as a BFGS-based realization of the generalized quadratic gradient principle.

We propose a unified framework that characterizes curvature-aware gradient optimization through positive-definite curvature transformations (GQG).

\subsection{Chiang's Quadratic Gradient}

Following the fixed-Hessian approach~\cite{bohning1988monotonicity}, where a simplified diagonal approximation is employed to reduce the computational complexity of second-order optimization, Bonte et al.~\cite{bonte2018privacy} applied this strategy to privacy-preserving machine learning. Building upon this line of research, Chiang~\cite{chiang2022privacy} introduced a more efficient gradient variant, referred to as the quadratic gradient, which incorporates curvature information into first-order optimization.

\subsubsection{ Original Quadratic Gradient} 
Given a differentiable scalar-valued function $F(\mathbf{x})$ with gradient $g$ and Hessian matrix $H$, Chiang~\cite{chiang2022privacy} introduced the quadratic gradient by constructing a suitable bound matrix of the Hessian. For maximization problems, a lower bound matrix $\bar{H}$ satisfying $\bar{H} \leq H$ is required, whereas for minimization problems, an upper bound matrix satisfying $H \leq \bar{H}$ is considered, where ``$\leq$'' denotes the Loewner ordering. The Hessian matrix $H$ itself satisfies these conditions and can therefore be directly used as the bound matrix. However, constructing a fixed bound matrix provides a more efficient alternative by avoiding repeated Hessian computations.

To derive the quadratic gradient, a diagonal scaling matrix $\bar{B}$ is first obtained from the bound matrix $\bar{H}$:

\begin{equation*}
  \begin{aligned}
   \bar B = 
\left[ \begin{array}{cccc}
  \frac{1}{ \epsilon + \sum_{i=0}^{d} | \bar h_{0i} | }   & 0  &  \ldots  & 0  \\
 0  &   \frac{1}{ \epsilon + \sum_{i=0}^{d} | \bar h_{1i} | }  &  \ldots  & 0  \\
 \vdots  & \vdots                & \ddots  & \vdots     \\
 0  &  0  &  \ldots  &   \frac{1}{ \epsilon + \sum_{i=0}^{d} | \bar h_{di} | }  \\
 \end{array}
 \right],
   \end{aligned}
\end{equation*}

where $\epsilon$ is a small positive constant introduced to prevent division by zero, and $\bar{h}_{ji}$ denotes the $(j,i)$-th element of $\bar{H}$. The quadratic gradient is then defined as:

\[
G=\bar{B}g.
\]

Since $\bar{B}$ is diagonal, the quadratic gradient $G$ has the same dimensionality as the original gradient $g$. Therefore, existing gradient-based optimization methods can be directly applied by replacing the conventional gradient with $G$, although a learning rate greater than $1$ is typically required. This formulation enables well-established first-order optimization techniques to incorporate second-order curvature information without explicitly computing the Hessian inverse.

\subsubsection{Simplified Quadratic Gradient} 
Although the original quadratic gradient exploits the information contained in every row of the bound matrix $\bar{H}$, the resulting computational cost becomes impractical for high-dimensional optimization problems, particularly in deep learning scenarios. To alleviate this issue, Chiang~\cite{chiang2026sqg} introduced Simplified Quadratic Gradient (SQG), which further reduces the complexity of curvature estimation by retaining only the diagonal elements of the bound matrix $\bar{H}$ (or the Hessian matrix $H$).

Specifically, the diagonal scaling matrix is simplified as:

\begin{equation*}
\bar B = \text{diag} \left( 
\frac{1}{\epsilon + |\bar{h}_{00}|}, 
\frac{1}{\epsilon + |\bar{h}_{11}|}, 
\ldots, 
\frac{1}{\epsilon + |\bar{h}_{dd}|} 
\right),
\end{equation*}

where $\bar{h}_{jj}$ denotes the $j$-th diagonal element of $\bar{H}$. The simplified quadratic gradient is subsequently computed as:

\[
G=\bar{B}g.
\]

By discarding the off-diagonal curvature terms, SQG reduces the construction of the quadratic gradient to element-wise operations, substantially decreasing memory consumption and computational overhead. Despite this simplification, SQG preserves optimization behavior comparable to the original quadratic gradient while providing improved scalability and compatibility with backpropagation-based large-scale stochastic optimization.

\subsubsection{Quasi-Quadratic Gradient} 
While the original quadratic gradient relies on a fixed or pre-defined Hessian bound matrix,
its static curvature approximation may fail to capture the dynamically evolving geometry of
complex optimization landscapes. To address this limitation, Chiang proposed the
Quasi-Quadratic Gradient (QQG) \citep{chiang2026qqg}, which extends the quadratic gradient
principle to the quasi-Newton framework by incorporating the adaptive curvature estimation
capability of BFGS. Unlike the original quadratic gradient that utilizes a fixed curvature matrix,
QQG directly leverages the inverse Hessian approximation maintained by BFGS.

Specifically, let $B_k$ denote the symmetric positive-definite (SPD) Hessian approximation
generated by the BFGS update at iteration $k$. The quasi-quadratic gradient is defined as:

\begin{equation*}
G_{\mathrm{qq}}^{(k)} = B_k^{-1} g_k,
\end{equation*}

where $g_k$ represents the current gradient and $B_k^{-1}$ provides an adaptive approximation
of the inverse Hessian. Since the BFGS update preserves the positive definiteness of $B_k$
under the curvature condition $s_k^T y_k > 0$, the resulting quasi-quadratic gradient
naturally incorporates local second-order geometric information while maintaining a reliable
optimization direction.

By replacing the static curvature approximation in the original quadratic gradient with the
dynamically updated BFGS inverse Hessian, QQG bridges the gap between first-order gradient
methods and quasi-Newton optimization. This formulation enables curvature-aware gradient
updates without explicitly computing the true Hessian matrix, providing an adaptive extension
of the quadratic gradient framework for complex optimization problems.

\subsubsection{Quadratic Gradient Algorithms}

The quadratic gradient (QG) framework provides a general mechanism for incorporating
curvature information into first-order optimization algorithms. Instead of directly
using the vanilla gradient, QG replaces the gradient update with a curvature-aware
direction, allowing existing first-order optimizers to benefit from second-order
information without explicitly computing the Hessian inverse.

Specifically, given the gradient $g=\nabla f(\mathbf{x})$, the quadratic gradient is
defined as:

\begin{equation}
    G = \bar{B}g,
\end{equation}

where $\bar{B}$ serves as an approximation of the inverse Hessian. Since the quadratic
gradient has the same dimensionality as the original gradient, it can be directly
integrated into existing gradient-based optimization algorithms.

Based on this property, Chiang~\cite{chiang2022privacy} demonstrated that the QG
framework can be applied to enhance various first-order optimization methods,
including Nesterov's Accelerated Gradient (NAG), AdaGrad, and Adam. By replacing the
vanilla gradient with the quadratic gradient, these enhanced optimizers incorporate
additional curvature information while preserving the original optimization structures.

For example, Enhanced NAG replaces the conventional gradient update:

\begin{equation}
    V_{t+1}=\boldsymbol{\beta}_t+\eta_t\nabla J(\boldsymbol{\beta}_t)
\end{equation}

with the quadratic-gradient update:

\begin{equation}
    V_{t+1}=\boldsymbol{\beta}_t+N_tG_t,
\end{equation}

where $G_t$ denotes the quadratic gradient and $N_t$ is the corresponding learning
rate scaling factor.

Similarly, Enhanced AdaGrad and Enhanced Adam substitute the original gradient
statistics with quadratic gradients. For AdaGrad, the parameter update becomes:

\begin{equation}
    \beta_i^{(t+1)}
    =
    \beta_i^{(t)}
    -
    \frac{N_t}
    {\epsilon+\sqrt{\sum_{k=1}^{t}(G_i^{(k)})^2}}
    G_i^{(t)}.
\end{equation}

For Adam, the first- and second-order moment estimates are computed using $G_t$
instead of the original gradient $g_t$, enabling adaptive optimization with additional
curvature information.

\subsection{Line Search Techniques}

Line search techniques play a complementary role in curvature-aware optimization
methods. While the proposed Generalized Quadratic Gradient (GQG) framework focuses on
constructing an improved optimization direction through a positive-definite
curvature transformation, line search methods are responsible for determining an
appropriate step length along this direction.

Therefore, GQG does not replace existing line search strategies. Instead, it can be
seamlessly combined with them: GQG improves the quality of the search direction,
while line search guarantees a suitable update magnitude and enhances global
convergence properties.

Given the GQG search direction

\begin{equation}
    p_k=-P_kg_k ,
\end{equation}

where $P_k\succ0$, a line search procedure determines the step size $\alpha_k$ and
updates the parameters as:

\begin{equation}
    x_{k+1}=x_k+\alpha_k p_k .
\end{equation}

Different line search strategies can be integrated into the GQG framework depending
on the desired trade-off between computational cost and convergence guarantees.

\paragraph{1. Exact Line Search}

Exact line search determines the optimal step size by minimizing the objective
function along the search direction:

\begin{equation}
    \alpha_k
    =
    \arg\min_{\alpha>0}f(x_k+\alpha p_k).
\end{equation}

Although exact line search provides the theoretically optimal step size for quadratic
objectives, it is computationally expensive for general nonlinear optimization
problems due to the large number of objective function evaluations required.
Consequently, practical implementations usually adopt inexact line search methods.

\paragraph{2. Inexact Line Search and Armijo Condition}

In practical optimization, it is often unnecessary to obtain the exact minimizer
along the search direction. Instead, inexact line search methods seek a step size
that provides sufficient decrease.

The Armijo condition requires:

\begin{equation}
    f(x_k+\alpha p_k)
    \leq
    f(x_k)+c_1\alpha\nabla f(x_k)^Tp_k ,
\end{equation}

where $c_1\in(0,1)$ is a small constant (typically $10^{-4}$).

This condition ensures that the selected step size achieves an adequate reduction in
the objective function while avoiding excessively large updates.

\paragraph{3. Wolfe Conditions}

For quasi-Newton-based implementations of GQG, the Wolfe conditions are particularly
important because they help preserve the positive definiteness of curvature updates.

The Wolfe conditions combine the Armijo sufficient decrease condition with the
curvature condition:

\begin{equation}
    \nabla f(x_k+\alpha p_k)^Tp_k
    \geq
    c_2\nabla f(x_k)^Tp_k ,
\end{equation}

where $c_1<c_2<1$.

The curvature condition prevents the step size from being too small and encourages
the satisfaction of:

\begin{equation}
    s_k^Ty_k>0 ,
\end{equation}

which is required to maintain the symmetric positive definiteness of BFGS-type
curvature approximations.

A stricter variant, known as the Strong Wolfe condition, replaces the curvature
condition with:

\begin{equation}
    |\nabla f(x_k+\alpha p_k)^Tp_k|
    \leq
    c_2|\nabla f(x_k)^Tp_k|.
\end{equation}

\paragraph{4. Backtracking Line Search}

Backtracking line search provides a computationally efficient implementation of
inexact line search. It starts from an initial step size (typically $\alpha=1$ for
Newton-type methods) and repeatedly reduces the step size by a factor
$\rho\in(0,1)$ until the Armijo condition is satisfied.

Because GQG already provides a curvature-aware search direction, backtracking line
search can be viewed as a lightweight mechanism for controlling the update magnitude.
This combination allows GQG to maintain the advantages of second-order geometric
information while avoiding unnecessary instability caused by inappropriate step
sizes.

\paragraph{Line Search Techniques}

It is important to emphasize that line search techniques and GQG address different
aspects of optimization. GQG determines \textit{where} to move by exploiting
positive-definite curvature information, whereas line search determines \textit{how
far} to move along this direction.

Consequently, line search remains an essential complementary component rather than a
replacement for GQG. The integration of GQG with adaptive step-size strategies
provides a flexible optimization framework that combines curvature-aware directions
with reliable convergence control.

\section{Methodology}

Although the Quadratic Gradient (QG) framework has demonstrated the effectiveness of
incorporating curvature information into first-order optimization, existing
formulations are mainly restricted to specific curvature constructions, such as
fixed Hessian approximations or their diagonal simplifications.

In this work, we revisit the fundamental principle behind quadratic gradient methods
and generalize it beyond predefined Hessian structures. Instead of limiting the
curvature transformation to a particular approximation, we introduce a unified
framework based on positive-definite curvature operators.

The proposed \textbf{Generalized Quadratic Gradient (GQG)} framework defines the
optimization direction as:

\begin{equation}
    G_k=P_k g_k,
\end{equation}

where $g_k=\nabla f(x_k)$ denotes the gradient and $P_k$ is a symmetric
positive-definite curvature transformation matrix.

Different choices of $P_k$ correspond to different realizations of the quadratic
gradient principle. The original QG can be interpreted as using a fixed curvature
approximation, while diagonal approximations lead to Simplified Quadratic Gradient
(SQG). Furthermore, by adopting the inverse Hessian approximation generated by
quasi-Newton methods as $P_k$, the framework naturally extends to the
Quasi-Quadratic Gradient (QQG).

\subsection{Generalized Quadratic Gradient Framework}
\label{subsec:gqg}
Second-order optimization methods are fundamentally based on exploiting the
local geometry of the objective function. Given a twice-differentiable
objective function $f(\mathbf{x})$, Newton-type methods approximate the
objective landscape around the current point using a local quadratic model:

\begin{equation}
    f(\mathbf{x}+\Delta\mathbf{x})
    \approx
    f(\mathbf{x})
    +
    g^T\Delta\mathbf{x}
    +
    \frac{1}{2}
    \Delta\mathbf{x}^T H\Delta\mathbf{x},
\end{equation}

where $g=\nabla f(\mathbf{x})$ and
$H=\nabla^2 f(\mathbf{x})$ denote the gradient and Hessian matrix,
respectively.

The optimal Newton direction is obtained by solving the stationary condition
of the quadratic model:

\begin{equation}
    H\Delta\mathbf{x}=-g.
\end{equation}

However, the exact Hessian matrix is not necessarily positive definite,
particularly in non-convex optimization problems. An indefinite Hessian may
produce non-descent directions and unstable optimization trajectories.
Therefore, constructing a reliable positive-definite substitute of the
Hessian has become one of the central problems in modern Newton-type
optimization.

The fundamental principle behind many second-order optimization algorithms is
to replace the original Hessian with a positive-definite curvature matrix:

\begin{equation}
    P_k \approx H_k^{-1},
    \qquad
    P_k \succ 0,
\end{equation}

which defines a stable local metric for transforming the gradient. From this
perspective, various Newton-type and quasi-Newton methods can be interpreted
as different strategies for constructing positive-definite curvature
operators.

This principle can be viewed as a variational interpretation of curvature
approximation. Instead of directly using the possibly indefinite Hessian,
optimization algorithms seek a positive-definite matrix that preserves the
essential geometric information of the local quadratic model while ensuring
a stable optimization direction. Typical constructions include fixed Hessian
approximations, diagonal approximations, quasi-Newton updates, Gauss-Newton
matrices, Fisher information matrices, and regularized Hessian modifications.

Based on this observation, we introduce the \textbf{Generalized Quadratic
Gradient (GQG)} framework, which abstracts the common principle behind these
methods. Given any positive-definite curvature operator $P_k$, the generalized
quadratic gradient is defined as:

\begin{equation}
    G_k=P_kg_k ,
\end{equation}

where $g_k$ is the first-order gradient. The matrix $P_k$ acts as a
curvature-aware transformation that reshapes the gradient according to the
local geometry of the objective function.

Unlike the original Quadratic Gradient framework, which relies on specific
fixed or diagonal Hessian approximations, GQG provides a general formulation
where different positive-definite curvature constructions lead to different
optimization algorithms. Therefore, the essential problem of designing a GQG
algorithm reduces to constructing an effective positive-definite curvature
matrix $P_k$.
\subsubsection{Motivation}
\label{subsec:motivation}

The conceptual foundation of the Quadratic Gradient (QG) framework originates from
the fixed-Hessian formulation of Newton-type optimization. By incorporating
second-order curvature information into gradient-based updates, QG provides a
computationally efficient approximation to Newton's method while preserving the
simplicity of first-order optimization algorithms.

However, existing QG formulations are mainly restricted to specific Hessian
approximations, such as fixed Hessian matrices or their diagonal simplifications.
Although these approximations reduce computational complexity, they may fail to
capture the evolving local geometry of modern high-dimensional and non-convex
optimization landscapes.

Meanwhile, the Quasi-Newton family, particularly the BFGS algorithm, provides an
adaptive mechanism for estimating curvature information through iterative secant
updates. Unlike fixed Hessian approaches, BFGS dynamically refines its Hessian
approximation using historical gradient information while maintaining a
symmetric positive-definite (SPD) structure.

These observations motivate us to revisit the fundamental principle behind quadratic
gradient methods. Rather than restricting the curvature transformation to a
particular Hessian approximation, we generalize the quadratic gradient concept by
allowing any valid positive-definite curvature operator to transform the gradient.

Specifically, we propose the \textbf{Generalized Quadratic Gradient (GQG)} framework:

\begin{equation}
    G_{gq}^{(k)}=P_k g_k,
\end{equation}

where $P_k$ represents a positive-definite curvature transformation matrix.

Different choices of $P_k$ recover different quadratic gradient variants. For
example, fixed Hessian approximations correspond to the original QG framework,
diagonal curvature approximations lead to Simplified Quadratic Gradient (SQG), and
the BFGS inverse Hessian approximation results in the Quasi-Quadratic Gradient
(QQG).

\subsubsection{Observation}
\label{subsec:observation}

A fundamental observation behind the Generalized Quadratic Gradient framework is that
the essential requirement for a curvature transformation is not a specific Hessian
construction, but the preservation of positive definiteness.

For a symmetric positive-definite matrix $P_k$, the transformed gradient direction
satisfies:

\begin{equation}
    g_k^T P_k g_k > 0 ,
\end{equation}

for any non-zero gradient $g_k$. Therefore, for minimization problems, the update
direction

\begin{equation}
    -P_k g_k
\end{equation}

is guaranteed to be a descent direction with respect to the local quadratic model.

This property provides the theoretical foundation of GQG and distinguishes it from
arbitrary gradient scaling methods. The curvature matrix $P_k$ does not need to be
the exact inverse Hessian; instead, it only needs to provide a meaningful
positive-definite approximation of the local geometry.

\paragraph{BFGS as a Positive-Definite Instance of GQG}

Among various curvature approximation techniques, BFGS provides a natural
instantiation of the GQG framework. Under the standard assumptions that the initial
matrix $B_0$ is symmetric positive definite and the curvature condition

\begin{equation}
    s_k^T y_k >0
\end{equation}

is satisfied, the BFGS update preserves the positive definiteness of the Hessian
approximation throughout the optimization process.

Consequently, the inverse Hessian approximation generated by BFGS:

\begin{equation}
    P_k=B_k^{-1},
\end{equation}

naturally satisfies the requirement of the GQG framework and leads to the
Quasi-Quadratic Gradient:

\begin{equation}
    G_{qq}^{(k)}
    =
    B_k^{-1}g_k .
\end{equation}

Therefore, QQG should be viewed as a BFGS-based realization of GQG rather than an
independent optimization paradigm.

\paragraph{Remark}

The positive definiteness requirement in GQG is independent of the optimization
objective. For minimization problems, a positive-definite curvature transformation
ensures a descent direction through $-P_kg_k$, while for maximization problems the
sign of the update is reversed.

This objective-independent formulation enables GQG to provide a unified perspective
for integrating various Hessian approximations and quasi-Newton techniques into
gradient-based optimization.

\subsubsection{Evolution of the Proposed Approach}
\label{subsec:evolution}

The development of the proposed Generalized Quadratic Gradient (GQG) framework
originated from a progressive investigation of how curvature information can be
incorporated into gradient-based optimization. Starting from the original Quadratic
Gradient (QG), we gradually identified that the essential principle is not restricted
to a specific Hessian approximation, but rather lies in the construction of an
appropriate positive-definite curvature transformation.

\paragraph{Initial Exploration: Extending Diagonal Quadratic Gradient}

The original Quadratic Gradient framework employs a curvature-aware gradient obtained
from a fixed Hessian-related approximation. Inspired by this formulation, our initial
attempts investigated whether more advanced curvature information, such as that
provided by quasi-Newton methods, could be incorporated into the existing diagonal
construction.

However, directly combining BFGS curvature estimation with the diagonal quadratic
gradient formulation resulted in unsatisfactory optimization behavior. This
observation suggested that the limitation was not necessarily caused by insufficient
curvature information, but rather by the restrictive diagonal approximation itself.
The diagonal constraint discards important correlations between optimization
dimensions and prevents the full utilization of the curvature structure captured by
quasi-Newton methods.

The diagonal restriction was an artificial limitation.

\paragraph{Beyond Diagonal Approximation: Positive-Definite Curvature Transformations}

The next stage of development explored replacing the diagonal scaling matrix with a
general positive-definite matrix. From the perspective of Newton-type optimization,
a full Hessian approximation naturally provides richer geometric information than
diagonal approximations.

However, directly adopting full Hessian matrices introduces two major challenges:

\begin{enumerate}
    \item \textbf{Computational Complexity:}
    Constructing and applying a full Hessian inverse is computationally expensive,
    especially for high-dimensional machine learning problems.

    \item \textbf{Curvature Construction:}
    Designing an efficient positive-definite approximation that provides accurate
    local geometry while maintaining numerical stability remains a fundamental
    challenge in second-order optimization.
\end{enumerate}

These observations motivated the investigation of a more general formulation in which
the curvature transformation is not restricted to a fixed Hessian approximation or a
specific diagonal structure.

\paragraph{From Quasi-Newton Approximation to Generalized Quadratic Gradient}

The key insight is that the quadratic gradient principle only requires a
positive-definite curvature transformation. Therefore, instead of designing a specific
replacement for the Hessian matrix, we formulate the generalized quadratic gradient as:

\begin{equation}
    G_k=P_kg_k,\qquad P_k\succ0 ,
\end{equation}

where $P_k$ represents a general curvature operator.

Under this formulation, different optimization techniques correspond to different
choices of $P_k$. The original Quadratic Gradient can be interpreted as using a fixed
curvature approximation, while the Quasi-Quadratic Gradient (QQG) corresponds to the
specific case where:

\begin{equation}
    P_k=H_k^{BFGS},
\end{equation}

with $H_k^{BFGS}$ denoting the inverse Hessian approximation generated by BFGS.

Therefore, QQG is not an independent optimization framework, but rather a
BFGS-based instance of the broader GQG framework.

The evolution of the proposed approach can be summarized as:

\begin{align*}
&\texttt{Newton-Raphson Method}
\xmapsto{}
\textit{Fixed Hessian Method}
\xmapsto{}
\textit{Simplified Fixed Hessian}
\\
&\xmapsto{}
\textit{Quadratic Gradient}
\xmapsto{}
\textit{Simplified Quadratic Gradient}
\\[3pt]
&\texttt{Newton-Raphson Method}
\xmapsto{}
\textit{Quasi-Newton Optimization}
\xmapsto{}
\textit{BFGS Inverse Hessian Approximation}
\\
&\xmapsto{}
\textit{Quasi-Quadratic Gradient}
\\[3pt]
&\textbf{Unified View:}
\\
&\texttt{Newton-Type Optimization}
\xmapsto{}
\textit{Positive-Definite Curvature Approximation}
\xmapsto{}
\textit{Generalized Quadratic Gradient}.
\end{align*}

Through this evolution, the proposed framework moves beyond designing individual
second-order optimization algorithms. Instead, GQG provides a unified perspective in
which various Hessian approximations and quasi-Newton methods can be interpreted as
different realizations of a common curvature-aware gradient principle.

GQG is not another optimizer; it is a unifying abstraction over curvature-aware optimization.


\subsubsection{Definition}

The original Quadratic Gradient (QG) framework introduces curvature information into
first-order optimization by replacing the vanilla gradient with a scaled gradient
obtained from a Hessian-related approximation. However, existing formulations mainly
rely on specific curvature constructions, such as fixed Hessian approximations,
diagonal simplifications, or quasi-Newton updates.

In this work, we generalize the quadratic gradient principle by considering a broader
class of positive-definite curvature transformations. Instead of restricting the
scaling matrix to a particular Hessian approximation, the proposed
\textbf{Generalized Quadratic Gradient (GQG)} defines the optimization direction through
a general positive-definite curvature matrix.

\begin{definition*}[\texttt{Generalized Quadratic Gradient}]
Let $g_k=\nabla F(\mathbf{x}_k)$ denote the gradient of the objective function at
iteration $k$. Let $P_k$ be a symmetric positive-definite (SPD) curvature matrix that
captures local geometric information of the optimization landscape. The Generalized
Quadratic Gradient is defined as:

\begin{equation*}
    G_{gq}^{(k)} = P_k g_k,
\end{equation*}

where $P_k\succ0$ represents a general curvature transformation.

Different choices of $P_k$ lead to different quadratic gradient variants. For example,
the original Quadratic Gradient corresponds to a fixed curvature approximation, while
the Quasi-Quadratic Gradient (QQG) can be obtained by selecting:

\begin{equation*}
    P_k = B_k^{-1},
\end{equation*}

where $B_k$ is the Hessian approximation maintained by the BFGS algorithm.
\end{definition*}

The generalized formulation reveals that the fundamental principle behind quadratic
gradient methods is not limited to a specific Hessian approximation, but rather the
construction of an appropriate positive-definite curvature operator. By transforming
the original gradient through $P_k$, GQG provides a unified framework that connects
first-order gradient methods with Newton-type and quasi-Newton optimization.

\subsubsection{Update Rules}

The Generalized Quadratic Gradient can be integrated into iterative optimization
algorithms through a simple gradient substitution. Given the parameter vector
$\beta$, the update rules are defined according to the optimization objective.

For maximization problems:

\begin{equation*}
    \beta_{k+1}
    =
    \beta_k+\eta_kG_{gq}^{(k)},
\end{equation*}

while for minimization problems:

\begin{equation*}
    \beta_{k+1}
    =
    \beta_k-\eta_kG_{gq}^{(k)},
\end{equation*}

where $\eta_k$ denotes the learning rate or step size.

By replacing the conventional gradient $g_k$ with the curvature-aware gradient
$G_{gq}^{(k)}$, GQG incorporates local second-order geometric information while
preserving the computational structure of first-order optimization methods. Therefore,
existing optimization algorithms can be naturally extended by replacing:

\begin{equation*}
    g_k \rightarrow G_{gq}^{(k)}.
\end{equation*}

\paragraph{Compatibility with Line Search Techniques}

Since GQG modifies only the search direction while preserving the general iterative
optimization framework, it is fully compatible with existing step-size selection
strategies. In particular, line search techniques such as backtracking line search and
Wolfe-condition-based methods can be directly integrated to determine $\eta_k$.

The curvature-aware direction provided by GQG and the adaptive step-size selection
provided by line search methods complement each other. The former improves the local
geometry of the optimization trajectory, while the latter ensures sufficient descent
and numerical stability.

Furthermore, GQG can be combined with momentum-based first-order methods, such as
Nesterov's Accelerated Gradient (NAG), as well as adaptive optimization algorithms
including AdaGrad and Adam. This compatibility enables GQG to serve as a general
curvature-enhancement mechanism for a wide range of optimization algorithms.

The Generalized Quadratic Gradient provides a unified framework for incorporating
positive-definite curvature information into gradient-based optimization. Depending
on the construction of the curvature matrix, different optimization algorithms can be
obtained as specific instances of GQG.
QQG is not a competing framework with GQG; rather, it is a particular realization of GQG obtained by choosing the BFGS inverse Hessian approximation as the positive-definite curvature operator.
Therefore, QQG can be interpreted as a BFGS-based implementation of GQG rather than an
independent optimization framework. The general GQG formulation allows the curvature
matrix to be constructed from any symmetric positive-definite approximation, while
QQG specifically exploits the adaptive inverse Hessian estimation provided by BFGS.
We identify a general principle behind quadratic-gradient methods and formulate a unified framework based on positive-definite curvature transformations. QQG is a natural BFGS-based realization of this framework.

\subsubsection{Variational Quadratic Model Principle}

The fundamental idea behind the Generalized Quadratic Gradient (GQG)
framework originates from the quadratic approximation underlying
Newton-type optimization.

At each iteration, second-order optimization methods construct a local
quadratic model:

\begin{equation}
m_k(d)
=
f(x_k)
+
g_k^T d
+
\frac{1}{2}d^T H_k d ,
\end{equation}

where $H_k$ represents the local curvature of the objective function.
However, the exact Hessian $\nabla^2 f(x_k)$ is often unavailable,
computationally expensive, or indefinite in non-convex regions.

Therefore, the central problem of curvature-aware optimization is to find
a positive-definite surrogate matrix:

\begin{equation}
\widetilde{H}_k \in \mathbb{S}_{++}^{n}
\end{equation}

that preserves essential geometric information while guaranteeing a stable
quadratic model.

We formulate this process as a variational curvature transformation:

\begin{equation}
\widetilde{H}_k
=
\mathcal{T}_{\mathrm{SPD}}
(\mathcal{I}_k),
\end{equation}

where $\mathcal{I}_k$ denotes available curvature information, including
Hessian approximations, gradient statistics, secant information, or
Jacobian structures, and $\mathcal{T}_{\mathrm{SPD}}$ maps this information
into the symmetric positive-definite space.

The resulting generalized quadratic gradient is defined as:

\begin{equation}
G_{\mathrm{GQG}}^{(k)}
=
\widetilde{H}_k^{-1}g_k .
\end{equation}

Therefore, GQG is not a specific optimization algorithm, but a universal
curvature transformation layer that converts arbitrary positive-definite
curvature approximations into Newton-like gradient directions.

\subsubsection{Positive-Definite Curvature Approximation}

The effectiveness of the GQG framework relies on constructing an appropriate
positive-definite curvature transformation. Unlike classical Newton methods, which
require explicitly computing and inverting the Hessian matrix at every iteration,
GQG allows flexible curvature approximations obtained from different optimization
techniques.

A desirable curvature transformation $P_k$ should satisfy:

\begin{equation}
    P_k\succ0,
\end{equation}

which guarantees that the update direction

\begin{equation}
    -P_k g_k
\end{equation}

forms a valid descent direction for minimization problems.

Among various curvature approximation techniques, quasi-Newton methods provide a
natural mechanism for constructing such positive-definite matrices through iterative
gradient-based updates.

\subsubsection{Relationship with Newton-Type Optimization}

The GQG framework can be viewed as a generalization of Newton-type optimization.
The classical Newton update is defined as

\begin{equation}
    \mathbf{x}_{k+1}
    =
    \mathbf{x}_k
    -
    H_k^{-1}g_k ,
\end{equation}

where $H_k=\nabla^2f(\mathbf{x}_k)$ is the Hessian matrix.

However, in non-convex optimization, the Hessian may be indefinite, making
$H_k^{-1}$ unreliable or even causing the Newton direction to become an ascent
direction. Therefore, practical second-order optimization methods replace the
original Hessian with a positive-definite approximation:

\begin{equation}
    P_k\approx H_k^{-1},
\end{equation}

leading to the generalized update:

\begin{equation}
    \mathbf{x}_{k+1}
    =
    \mathbf{x}_k-\eta_kG_k ,
\end{equation}

where $\eta_k$ is the learning rate.

Therefore, GQG provides a unified formulation covering Newton, quasi-Newton,
and adaptive optimization methods through different choices of $P_k$.

\subsubsection{Connection to Existing Quadratic Gradient Methods}

The original Quadratic Gradient can be interpreted as a special case of GQG.
Specifically, when a fixed positive-definite approximation $\bar{H}$ is used,
the curvature transformation becomes

\begin{equation}
    P_k=\bar{H}^{-1},
\end{equation}

and the generalized quadratic gradient reduces to

\begin{equation}
    G_k=\bar{H}^{-1}g_k .
\end{equation}

The Simplified Quadratic Gradient further restricts $P_k$ to a diagonal
structure:

\begin{equation}
    P_k=
    \operatorname{diag}(p_1,\ldots,p_d).
\end{equation}

Although these approaches differ in their curvature construction strategies,
they share the same fundamental principle:

\begin{equation}
\boxed{
    \text{Gradient direction}
    \times
    \text{Positive-definite curvature transformation}
    =
    \text{Quadratic Gradient}
}
\end{equation}

\subsubsection{Stationary Quadratic Model Interpretation}

The motivation behind GQG originates from the local quadratic approximation of
the objective function:

\begin{equation}
f(\mathbf{x}+\Delta\mathbf{x})
\approx
f(\mathbf{x})
+
g^T\Delta\mathbf{x}
+
\frac12
\Delta\mathbf{x}^TH\Delta\mathbf{x}.
\end{equation}

The stationary point of this local quadratic model satisfies

\begin{equation}
    H\Delta\mathbf{x}=-g.
\end{equation}

When the exact Hessian is unavailable or unreliable, GQG replaces the inverse
Hessian by a positive-definite curvature surrogate:

\begin{equation}
    \Delta\mathbf{x}
    =
    -P_kg.
\end{equation}

Consequently, the essential objective of GQG is not to reproduce the exact
Hessian, but to construct a reliable positive-definite approximation that
captures meaningful curvature information while maintaining numerical stability.

\subsubsection{Generalization Principle}

The proposed framework reveals that existing quadratic-gradient algorithms can
be interpreted as different realizations of the same principle:

\begin{equation}
    \textbf{Curvature Construction}
    \rightarrow
    P_k\succ0
    \rightarrow
    G_k=P_kg_k .
\end{equation}

Under this perspective:

\begin{itemize}
    \item Original Quadratic Gradient corresponds to fixed Hessian-based
    curvature construction;
    \item Simplified Quadratic Gradient corresponds to diagonal curvature
    approximation;
    \item Quasi-Quadratic Gradient corresponds to BFGS-based dynamic curvature
    approximation;
    \item GQG-Adam and GQG-AdaGrad correspond to combining curvature-aware
    gradients with adaptive optimization mechanisms.
\end{itemize}

Therefore, the Generalized Quadratic Gradient framework provides a unified view
of curvature-aware optimization, where the essential design problem is shifted
from a specific optimizer formulation to the construction of a reliable
positive-definite curvature transformation.

\subsection{Positive-Definite Curvature Construction}
\label{subsec:positive_definite_curvature}

A fundamental challenge in second-order optimization is that the Hessian matrix
$\nabla^2 f(\mathbf{x})$ is generally indefinite, especially in non-convex
optimization landscapes. Although the Newton direction provides fast local
convergence near optimal solutions, the presence of negative eigenvalues may
result in non-descent directions and numerical instability. Therefore, a central
problem in modern optimization is how to construct a reliable positive-definite
curvature matrix that preserves useful second-order information while ensuring
stable optimization dynamics.

In the proposed Generalized Quadratic Gradient (GQG) framework, we abstract this
problem by considering a general positive-definite curvature transformation:
\begin{equation}
    G_k = P_k g_k,
\end{equation}
where $g_k$ denotes the first-order gradient and $P_k$ represents a
positive-definite curvature matrix or preconditioner:
\begin{equation}
    P_k \succ 0.
\end{equation}

Under this formulation, existing Hessian approximation techniques can be
interpreted as different strategies for constructing $P_k$. The proposed GQG
framework therefore provides a unified perspective that connects Newton-type
methods, quasi-Newton methods, adaptive gradient methods, and modern
curvature-aware optimization algorithms.

\subsubsection{Quasi-Newton-Based Curvature Approximation}

Quasi-Newton methods represent one of the most successful approaches for
constructing positive-definite Hessian approximations without explicitly
computing second-order derivatives. Instead of directly evaluating
$\nabla^2 f(\mathbf{x})$, quasi-Newton algorithms iteratively update a curvature
matrix using only gradient differences.

Among various quasi-Newton methods, BFGS is particularly important due to its
ability to preserve the symmetric positive-definite property. Given
\begin{equation}
    s_k=x_{k+1}-x_k,
    \qquad
    y_k=g_{k+1}-g_k,
\end{equation}
the BFGS update constructs

\begin{equation}
B_{k+1}
=
B_k
+
\frac{y_ky_k^T}{y_k^Ts_k}
-
\frac{B_ks_ks_k^TB_k}{s_k^TB_ks_k}.
\end{equation}

When the initial matrix $B_0$ is positive definite and the curvature condition
$s_k^Ty_k>0$ is satisfied, the sequence $\{B_k\}$ remains symmetric positive
definite. Therefore, the inverse Hessian approximation
$P_k=B_k^{-1}$ naturally satisfies the requirement of the GQG framework.

The resulting update direction,

\begin{equation}
    G_k=B_k^{-1}g_k ,
\end{equation}

corresponds to the Quasi-Quadratic Gradient (QQG), which can be viewed as a
specific BFGS-based realization of the generalized quadratic gradient principle.

\subsubsection{Gauss-Newton and Fisher-Based Curvature Approximation}

Another important family of positive-definite curvature construction methods is
based on gradient outer products. These approaches replace the potentially
indefinite Hessian with naturally positive semi-definite matrices.

For nonlinear least-squares problems, the Gauss-Newton method approximates the
Hessian using

\begin{equation}
    G=J^TJ,
\end{equation}

where $J$ denotes the Jacobian matrix. Since $J^TJ$ is a Gram matrix, it is
guaranteed to satisfy

\begin{equation}
    J^TJ\succeq0.
\end{equation}

The Generalized Gauss-Newton (GGN) method extends this idea to broader classes
of machine learning objectives by replacing the Hessian with a positive
semi-definite approximation of the loss curvature.

Similarly, the Fisher Information Matrix (FIM) constructs curvature information
from the statistical geometry of the model:

\begin{equation}
    F=
    \mathbb{E}
    [
    gg^T
    ],
\end{equation}

which is always positive semi-definite. The empirical Fisher approximation,

\begin{equation}
    F\approx
    \frac1N
    \sum_{i=1}^{N}g_i g_i^T ,
\end{equation}

can be interpreted as a gradient-based curvature approximation, where the
curvature information is obtained directly from accumulated gradient statistics.

These methods demonstrate that positive-definite curvature can be constructed
without explicit Hessian evaluation, motivating the broader GQG perspective.

\subsubsection{Regularized Hessian Modification}

When the original Hessian contains negative eigenvalues, another common strategy
is to modify the Hessian directly to obtain a positive-definite approximation.

A representative example is the Levenberg--Marquardt (LM) method, which applies
a diagonal shift:

\begin{equation}
    G_\lambda = J^TJ+\lambda I .
\end{equation}

The regularization term $\lambda I$ guarantees positive definiteness while
maintaining the original curvature structure.

Similarly, damped Newton methods modify the Hessian as

\begin{equation}
    H_\lambda=H+\lambda I,
\end{equation}

where $\lambda$ is chosen to compensate for negative curvature.

More sophisticated approaches, such as Modified Cholesky factorization, directly
decompose the Hessian:

\begin{equation}
    H\approx LDL^T,
\end{equation}

and modify unstable diagonal elements or negative pivots to construct an SPD
approximation.

Another spectral approach is eigenvalue modification. Given

\begin{equation}
    H=U\Lambda U^T,
\end{equation}

negative eigenvalues can be clipped:

\begin{equation}
    \Lambda_i
    \leftarrow
    \max(\Lambda_i,\epsilon),
\end{equation}

producing a positive-definite matrix while preserving the original eigenspace.

\subsection{Generalized Quadratic Gradient Algorithms}

The integration of Generalized Quadratic Gradient (GQG) into first-order optimization
methods can be achieved by replacing the conventional gradient with a curvature-aware
gradient direction. Given the gradient $g_t$ and a positive-definite curvature matrix
$B_t$, the generalized quadratic gradient is defined as:

\begin{equation}
    G_t = B_t g_t,
\end{equation}

where $B_t \succ 0$ represents a general curvature approximation. Different choices of
$B_t$ lead to different instances of GQG, including fixed Hessian approximations,
simplified diagonal approximations, and quasi-Newton approximations such as BFGS.

\begin{enumerate}

\item \textbf{GQG-NAG:}

The integration of GQG with momentum-based optimization requires careful consideration
of the interaction between momentum accumulation and curvature-aware updates. The
standard Nesterov accelerated gradient can be modified by replacing the vanilla
gradient with $G_t$:

\begin{align*}
    V_{t+1} &= \boldsymbol{\beta}_t+\eta_t G_t,\\
    \boldsymbol{\beta}_{t+1}
    &=
    (1-\gamma_t)V_{t+1}+\gamma_t V_t,
\end{align*}

where $\eta_t$ controls the effective step size. Since the scale of $G_t$ depends on
the selected curvature matrix, adaptive step-size strategies, such as warm-up
scheduling or line-search methods, can be employed to improve stability.

\item \textbf{GQG-AdaGrad:}

AdaGrad can be extended by accumulating the generalized quadratic gradients instead of
the original gradients. The update rule becomes:

\begin{equation*}
    \beta_i^{(t+1)}
    =
    \beta_i^{(t)}
    -
    \frac{\eta}
    {\epsilon+\sqrt{\sum_{k=1}^{t}(G_i^{(k)})^2}}
    G_i^{(t)},
\end{equation*}

where $G_i^{(t)}$ denotes the $i$-th component of the generalized quadratic gradient.
This formulation preserves the adaptive learning-rate mechanism of AdaGrad while
incorporating curvature information through $B_t$.

\item \textbf{GQG-Adam:}

Following the same principle, GQG can be integrated into Adam by replacing the
original gradient in the first- and second-moment estimations:

\begin{align*}
    m_t &=
    \beta_1m_{t-1}
    +(1-\beta_1)G_t,\\
    v_t &=
    \beta_2v_{t-1}
    +(1-\beta_2)G_t^2.
\end{align*}

The remaining bias correction and parameter update procedures remain unchanged.
Therefore, GQG-Adam inherits Adam's adaptive optimization mechanism while utilizing
curvature-aware update directions.
\end{enumerate}

To evaluate the practical behavior of GQG-enhanced optimizers, we consider different
curvature matrix constructions and corresponding hyperparameter configurations. Since
the magnitude of the generalized quadratic gradient depends on the selected curvature
approximation, the learning rate may require adjustment compared with conventional
first-order optimization.

For example, when $B_t$ provides an approximate inverse-Hessian scaling, larger
learning rates may become applicable due to the improved conditioning of the update
direction. The influence of curvature construction and step-size selection is
investigated experimentally in the following section.

\subsubsection{Example: GQG-Adam}

The Adam optimizer~\cite{kingma2014adam} is governed by three primary hyper-parameters,
each controlling a different aspect of adaptive gradient optimization:

\begin{itemize}
    \item \textbf{Learning Rate ($\alpha$):} 
    The learning rate determines the magnitude of parameter updates and directly affects
    the convergence speed and optimization stability. In Adam, $\alpha$ controls the
    overall scale of the adaptive update.

    \item \textbf{First Moment Decay ($\beta_1$):}
    The parameter $\beta_1$ controls the exponential moving average of gradients, which
    introduces momentum into the optimization process. A larger value of $\beta_1$
    provides stronger smoothing of stochastic gradient noise and encourages updates
    along consistent optimization directions.

    \item \textbf{Second Moment Decay ($\beta_2$):}
    The parameter $\beta_2$ controls the exponential moving average of squared gradients.
    This adaptive mechanism adjusts the update magnitude according to historical gradient
    statistics, allowing parameters with different gradient scales to be optimized with
    appropriate learning rates.
\end{itemize}

To demonstrate the general applicability of the proposed Generalized Quadratic Gradient
(GQG), we integrate GQG into the Adam optimizer, resulting in a curvature-aware variant
denoted as \textit{GQG-Adam}. The key modification is to replace the vanilla gradient
$g_t$ in Adam with the generalized quadratic gradient:

\begin{equation}
    G_t = B_t g_t,
\end{equation}

where $B_t$ represents a positive-definite curvature matrix. Depending on the selected
curvature construction strategy, $B_t$ can represent different Newton-type or
quasi-Newton approximations.

After this substitution, the fundamental roles of Adam's hyper-parameters remain
unchanged:

\begin{itemize}
    \item \textbf{Momentum Smoothing ($\beta_1$):}
    In GQG-Adam, $\beta_1$ controls the exponential moving average of the generalized
    quadratic gradients. It preserves the momentum mechanism of Adam while incorporating
    curvature-aware optimization directions.

    \item \textbf{Curvature-aware Adaptation ($\beta_2$):}
    The parameter $\beta_2$ accumulates the squared magnitude of the generalized
    quadratic gradients instead of the original gradients. Therefore, Adam's adaptive
    scaling mechanism operates on curvature-enhanced update directions while maintaining
    its parameter-wise learning rate adaptation.

    \item \textbf{Effective Step-size Control ($\alpha$):}
    The learning rate $\alpha$ continues to regulate the global update scale. Since GQG
    modifies the gradient direction through the curvature matrix rather than changing
    Adam's optimization structure, conventional learning-rate tuning strategies remain
    applicable.
\end{itemize}

\section{Experiments}

In this section, we evaluate the effectiveness of Generalized Quadratic Gradient (GQG)
under different positive-definite curvature constructions. The objective of these
experiments is not only to compare optimization performance with first-order methods,
but also to investigate whether the quadratic gradient principle can be generalized
beyond specific Hessian approximation strategies.

We consider both convex and non-convex benchmark functions with different curvature
structures. These experiments are designed to evaluate convergence behavior, robustness
to ill-conditioned landscapes, and the influence of different curvature matrix
constructions.

\subsection{Convex Benchmarks}

We first evaluate GQG on convex optimization problems to analyze its convergence
behavior under well-defined curvature structures. The following benchmark functions
are considered:

\begin{itemize}
    \item \textbf{Sphere Function:}
    \begin{equation}
        f(\mathbf{x})=\sum_{i=1}^{n}x_i^2,
    \end{equation}
    which provides a simple isotropic convex landscape.

    \item \textbf{Sum of Different Powers Function:}
    \begin{equation}
        f(\mathbf{x})=\sum_{i=1}^{n}|x_i|^{i+1},
    \end{equation}
    which introduces varying curvature across dimensions and evaluates the ability of
    curvature-aware methods to handle heterogeneous optimization landscapes.
\end{itemize}

\subsection{Non-Convex Benchmarks}

To evaluate the robustness of GQG in complex optimization landscapes, we further
consider several non-convex benchmark functions:

\begin{itemize}
    \item \textbf{Rosenbrock Function:}
    \begin{equation}
        f(\mathbf{x})
        =
        \sum_{i=1}^{n-1}
        [
        100(x_{i+1}-x_i^2)^2+(1-x_i)^2
        ],
    \end{equation}
    which contains a narrow curved valley and is widely used to evaluate optimization
    methods under ill-conditioned curvature.

    \item \textbf{Rastrigin Function:}
    \begin{equation}
        f(\mathbf{x})
        =
        10n+\sum_{i=1}^{n}
        [x_i^2-10\cos(2\pi x_i)],
    \end{equation}
    which introduces multiple local minima and tests the robustness of optimization
    methods in highly non-convex landscapes.
\end{itemize}

\subsection{Saddle Point Analysis}

Saddle points are challenging for gradient-based optimization because the gradient
magnitude can become small even when the current point is far from an optimum. To
investigate the behavior of curvature-aware updates around saddle regions, we consider
the Monkey Saddle function:

\begin{equation}
    f(x,y)=x^3-3xy^2.
\end{equation}

The Hessian of this function is indefinite around the saddle point, providing a useful
test case for analyzing the effect of positive-definite curvature approximations.

Unlike first-order methods that rely solely on gradient magnitude, GQG incorporates
local curvature information through a positive-definite matrix $B$:

\begin{equation}
    G=B\nabla f(\mathbf{x}).
\end{equation}

Therefore, different choices of $B$ lead to different transformations of the local
optimization geometry. This experiment investigates how various curvature
constructions influence optimization trajectories near saddle regions.

\subsection{Evaluation Protocol}

For each benchmark function, we compare different GQG variants constructed from
different positive-definite curvature matrices. The evaluation metrics include:

\begin{itemize}
    \item convergence rate measured by objective reduction over iterations;
    \item number of iterations required to reach a predefined tolerance;
    \item robustness under different initialization points.
\end{itemize}

The experimental results provide an empirical analysis of how curvature construction
strategies affect the performance of generalized quadratic gradient optimization.

\begin{figure}[htbp]
\centering
\captionsetup[subfigure]{justification=centering}

\subfloat[The iDASH dataset]{%
    \includegraphics[width=0.48\textwidth]{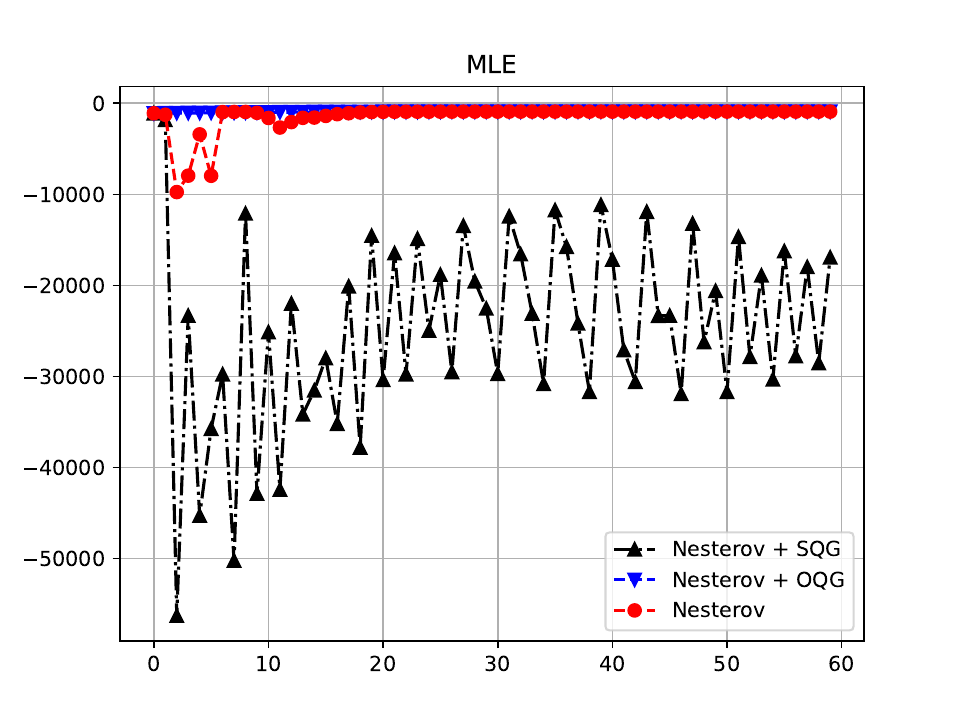}
    \label{fig:subfig01}
}
\hfill
\subfloat[The Edinburgh datasetn]{%
    \includegraphics[width=0.48\textwidth]{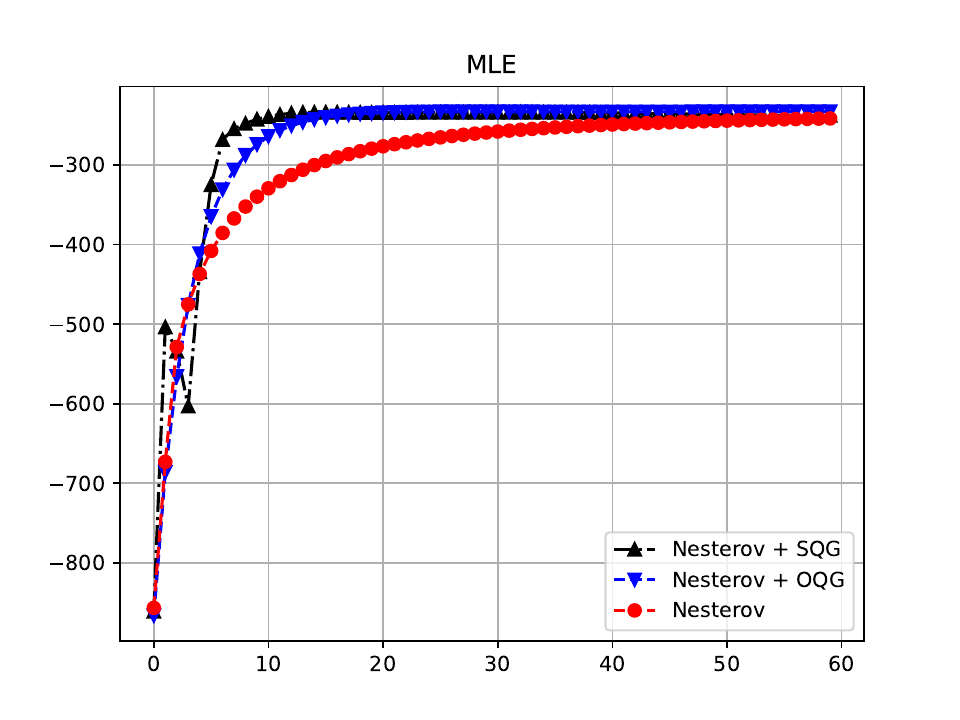}
    \label{fig:subfig02}
}

\vspace{1em} 

\subfloat[The lbw dataset]{%
    \includegraphics[width=0.48\textwidth]{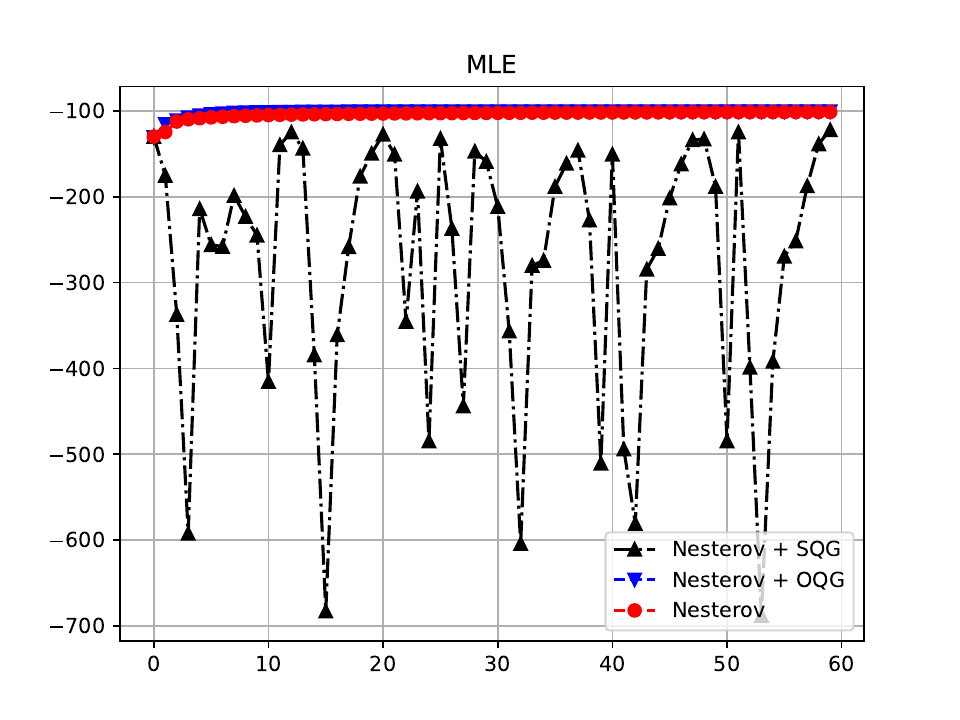}
    \label{fig:subfig03}
}
\hfill
\subfloat[The nhanes3 dataset]{%
    \includegraphics[width=0.48\textwidth]{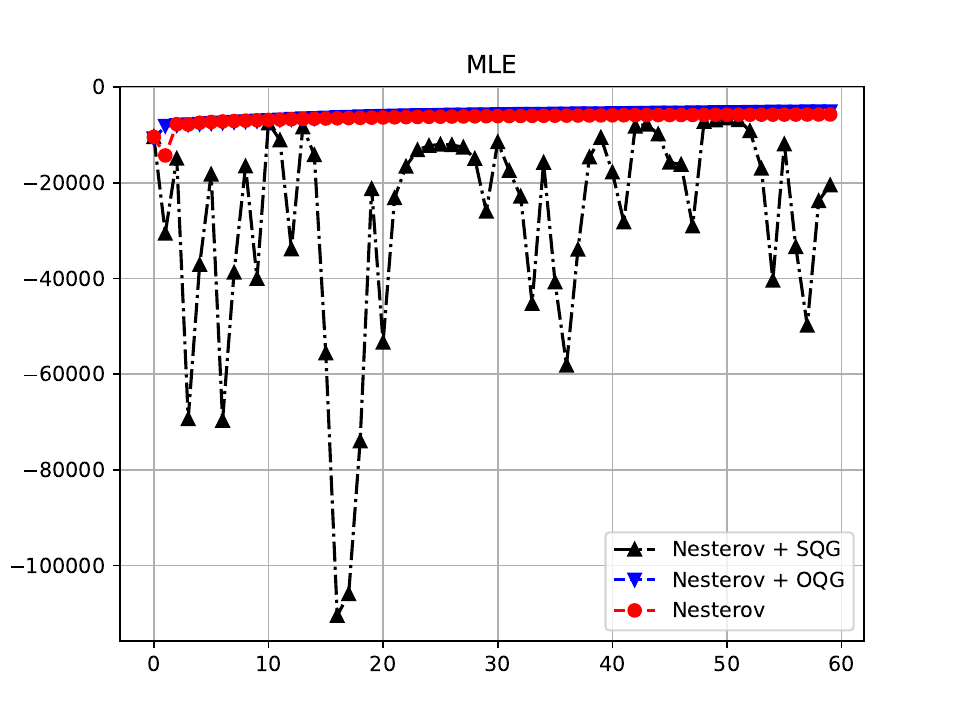}
    \label{fig:subfig04}
}

\vspace{1em} 

\subfloat[The pcs dataset]{%
    \includegraphics[width=0.48\textwidth]{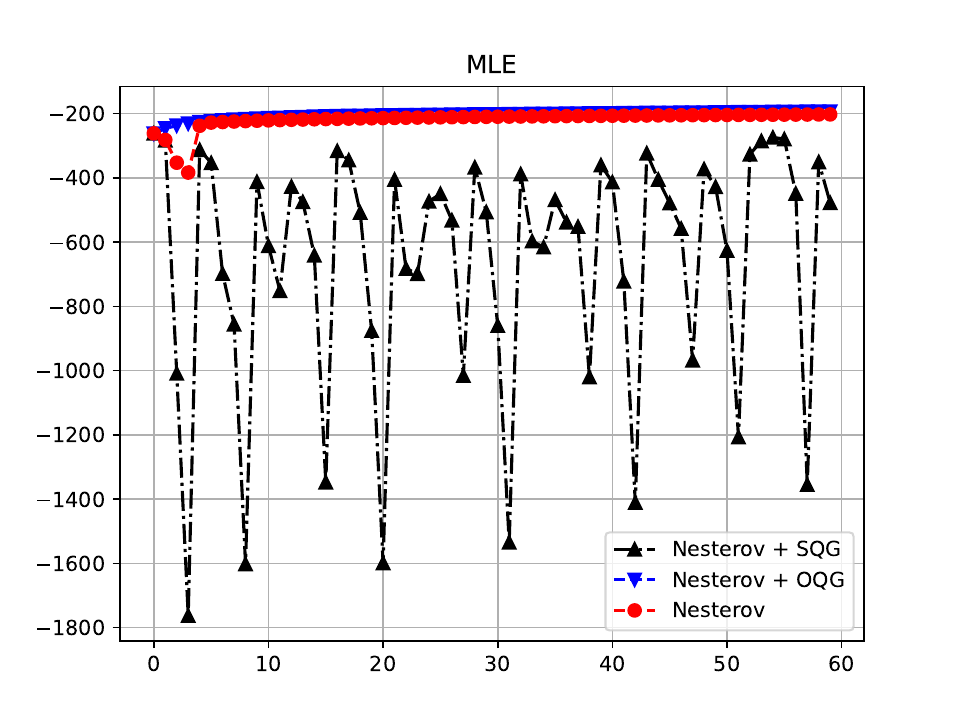}
    \label{fig:subfig03}
}
\hfill
\subfloat[The uis dataset]{%
    \includegraphics[width=0.48\textwidth]{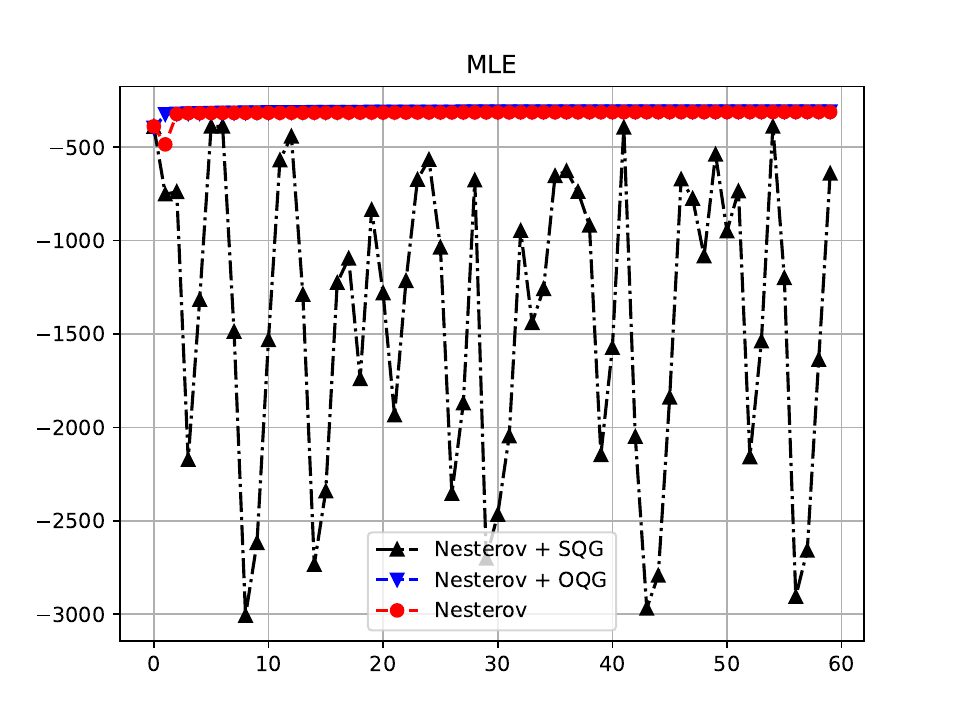}
    \label{fig:subfig04}
}

\vspace{1em} 

\subfloat[restructured MNIST dataset]{%
    \includegraphics[width=0.48\textwidth]{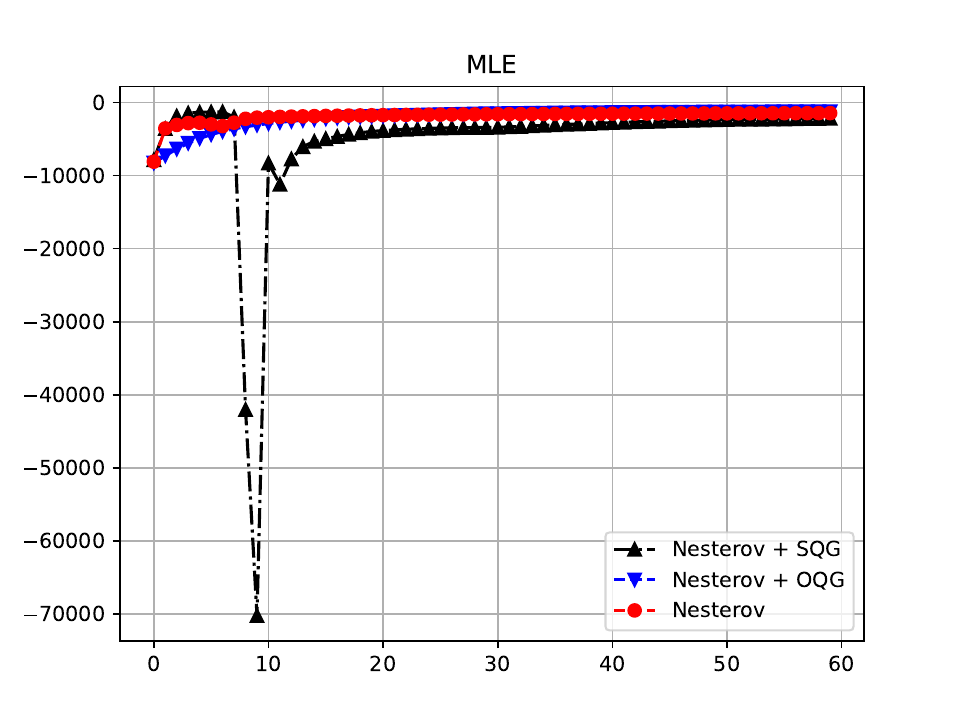}
    \label{fig:subfig03}
}
\hfill
\subfloat[The private financial dataset]{%
    \includegraphics[width=0.48\textwidth]{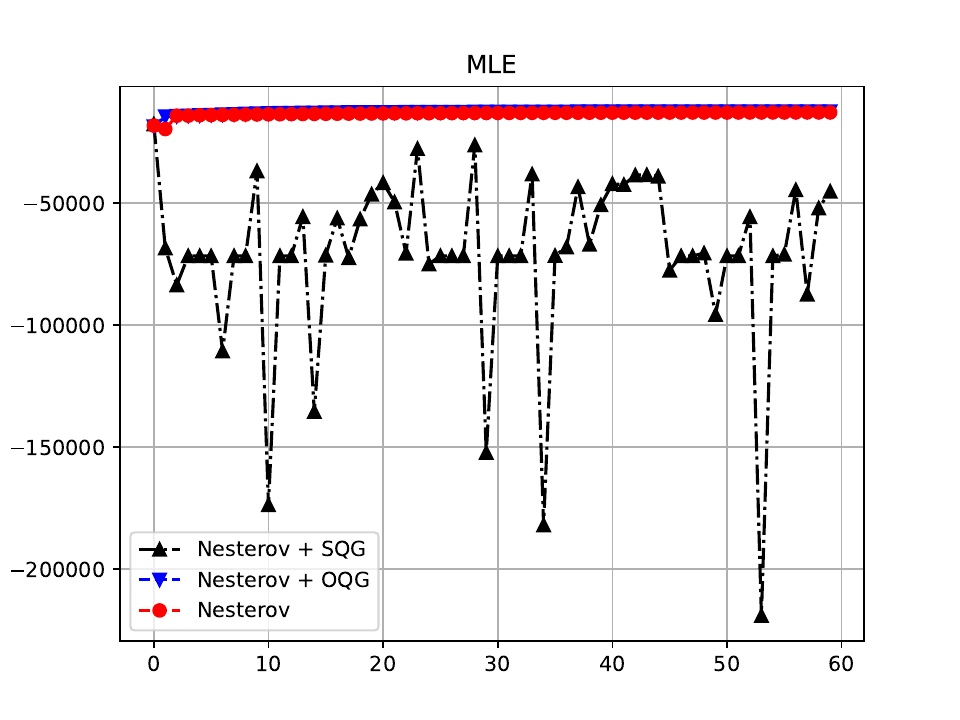}
    \label{fig:subfig04}
}

\caption{The training results of NAG + SQG vs. NAG + OQG vs. NAG in the clear domain.}
\label{fig0}
\end{figure}

\begin{figure}[htbp]
\centering
\captionsetup[subfigure]{justification=centering}

\subfloat[The iDASH dataset]{%
    \includegraphics[width=0.48\textwidth]{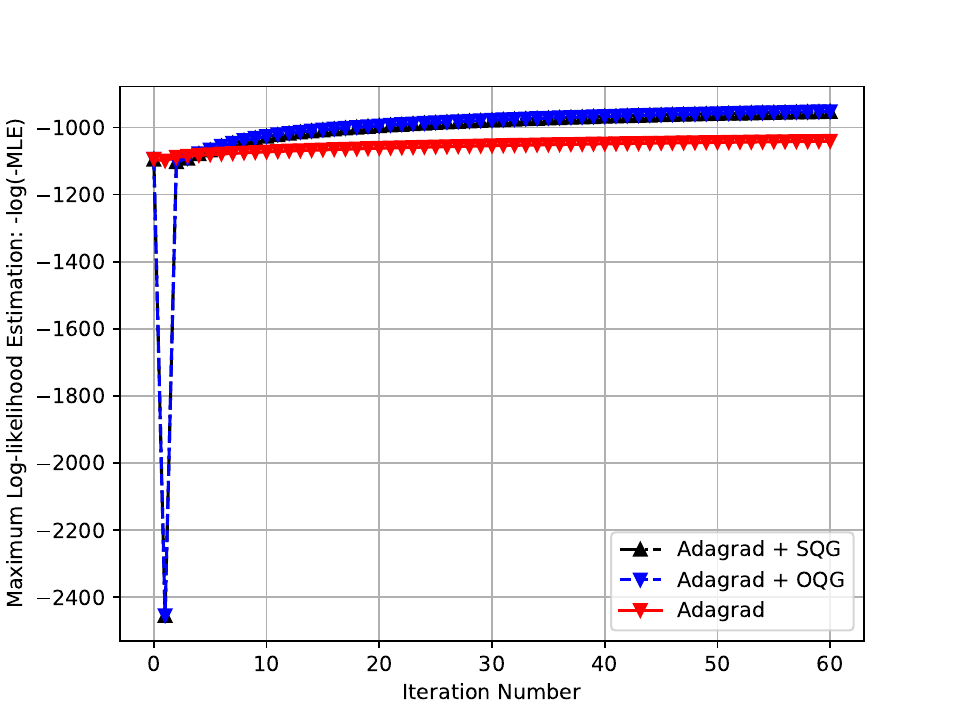}
    \label{fig:subfig01}
}
\hfill
\subfloat[The Edinburgh datasetn]{%
    \includegraphics[width=0.48\textwidth]{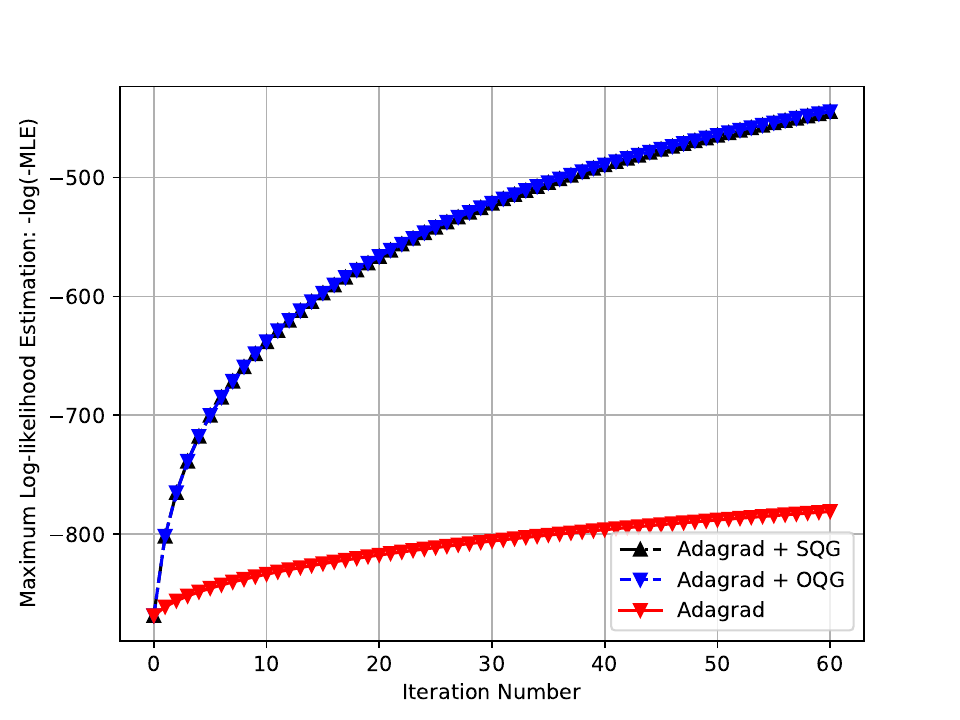}
    \label{fig:subfig02}
}

\vspace{1em} 

\subfloat[The lbw dataset]{%
    \includegraphics[width=0.48\textwidth]{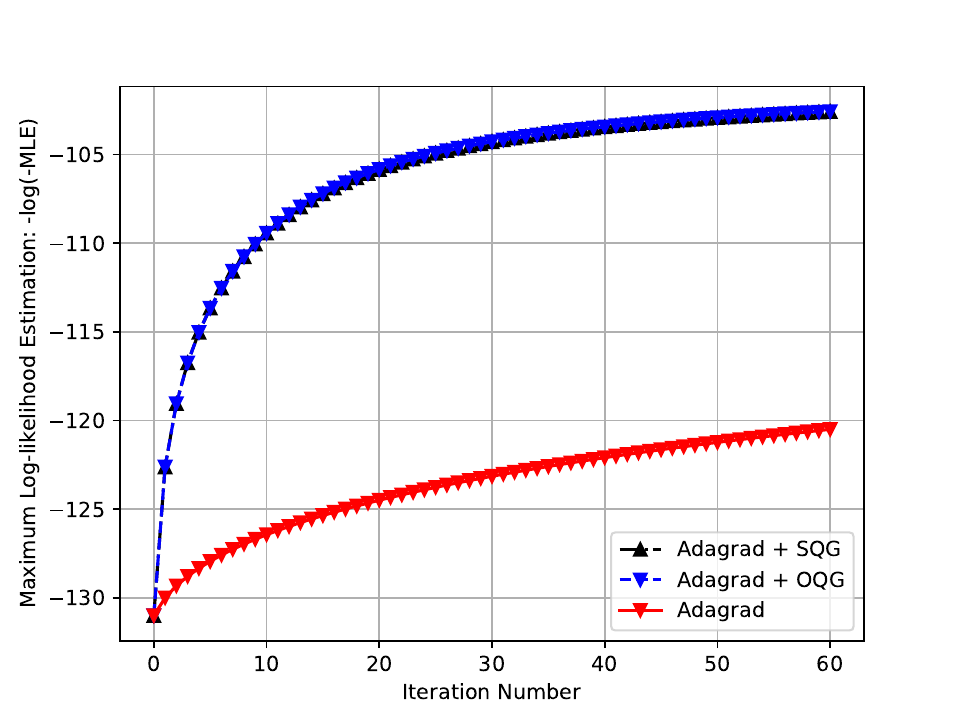}
    \label{fig:subfig03}
}
\hfill
\subfloat[The nhanes3 dataset]{%
    \includegraphics[width=0.48\textwidth]{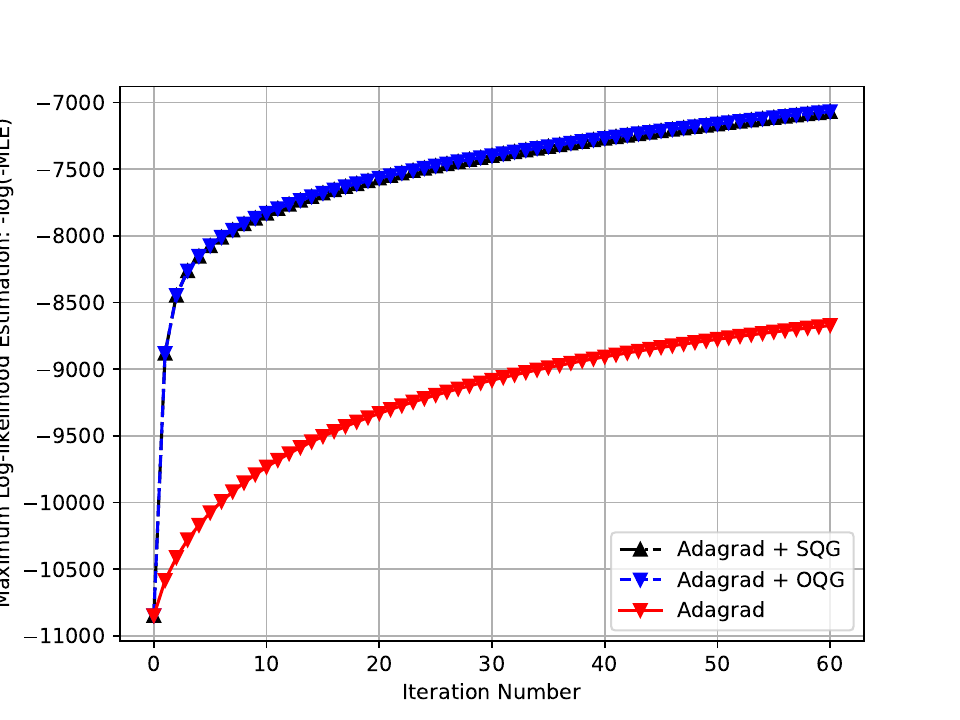}
    \label{fig:subfig04}
}

\vspace{1em} 

\subfloat[The pcs dataset]{%
    \includegraphics[width=0.48\textwidth]{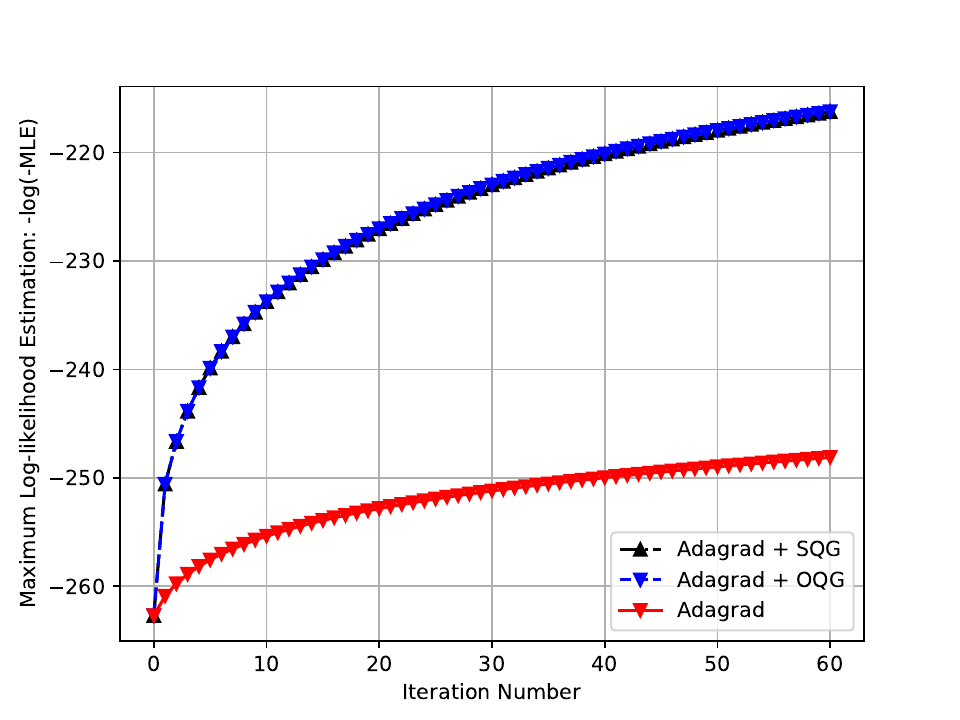}
    \label{fig:subfig03}
}
\hfill
\subfloat[The uis dataset]{%
    \includegraphics[width=0.48\textwidth]{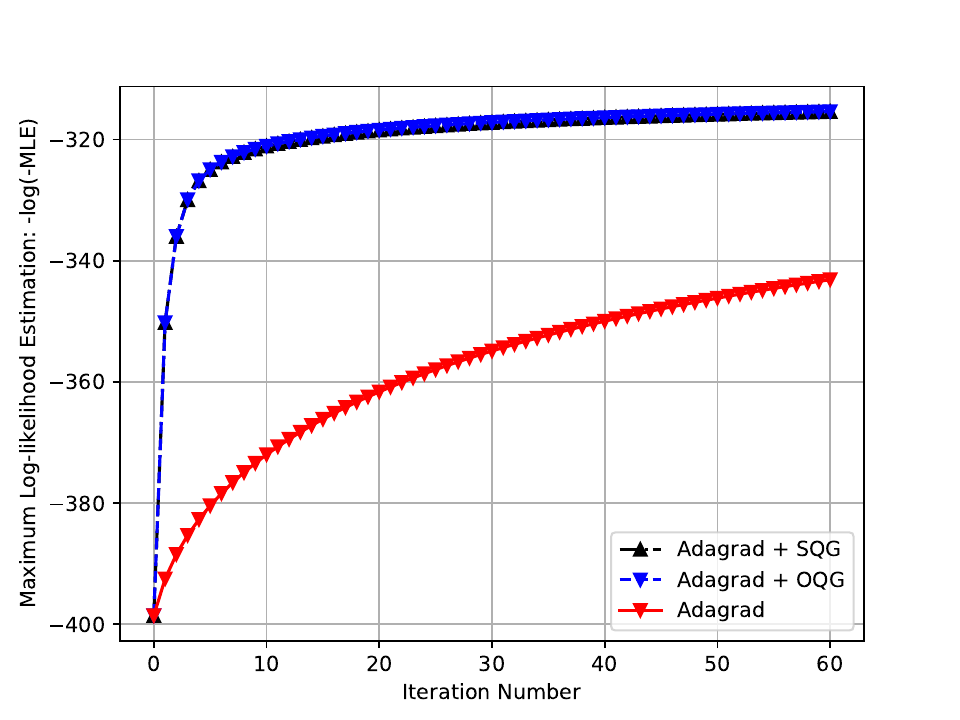}
    \label{fig:subfig04}
}

\vspace{1em} 

\subfloat[restructured MNIST dataset]{%
    \includegraphics[width=0.48\textwidth]{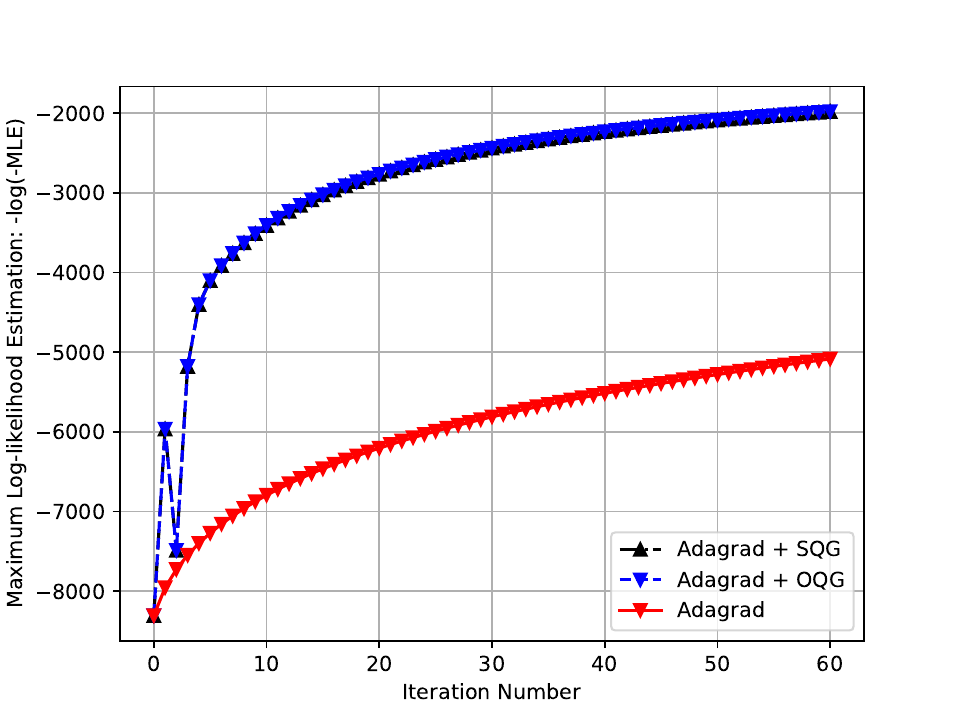}
    \label{fig:subfig03}
}
\hfill
\subfloat[The private financial dataset]{%
    \includegraphics[width=0.48\textwidth]{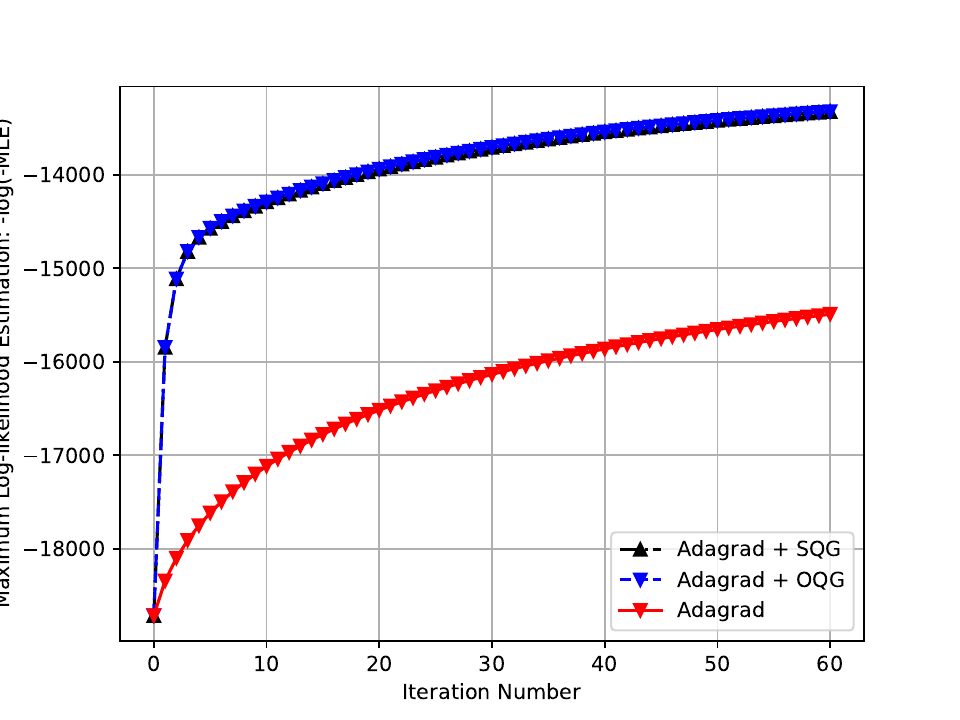}
    \label{fig:subfig04}
}

\caption{The training results of AdaGrad + SQG vs. AdaGrad + OQG vs. AdaGrad in the clear domain.}
\label{fig1}
\end{figure}

\begin{figure}[htbp]
\centering
\captionsetup[subfigure]{justification=centering}

\subfloat[The iDASH dataset]{%
    \includegraphics[width=0.48\textwidth]{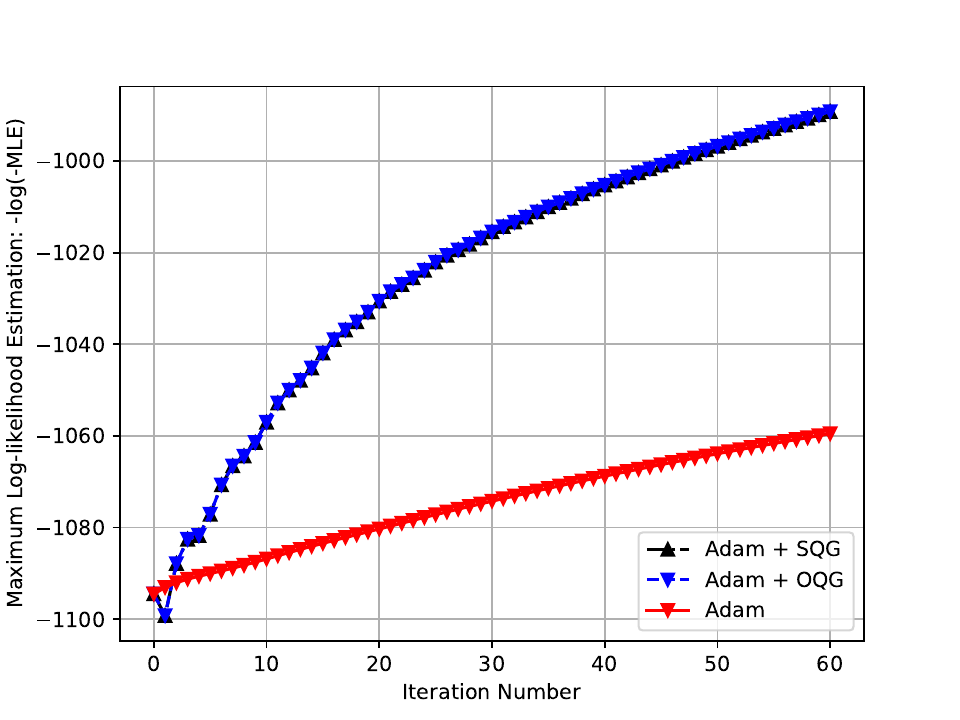}
    \label{fig:subfig01}
}
\hfill
\subfloat[The Edinburgh datasetn]{%
    \includegraphics[width=0.48\textwidth]{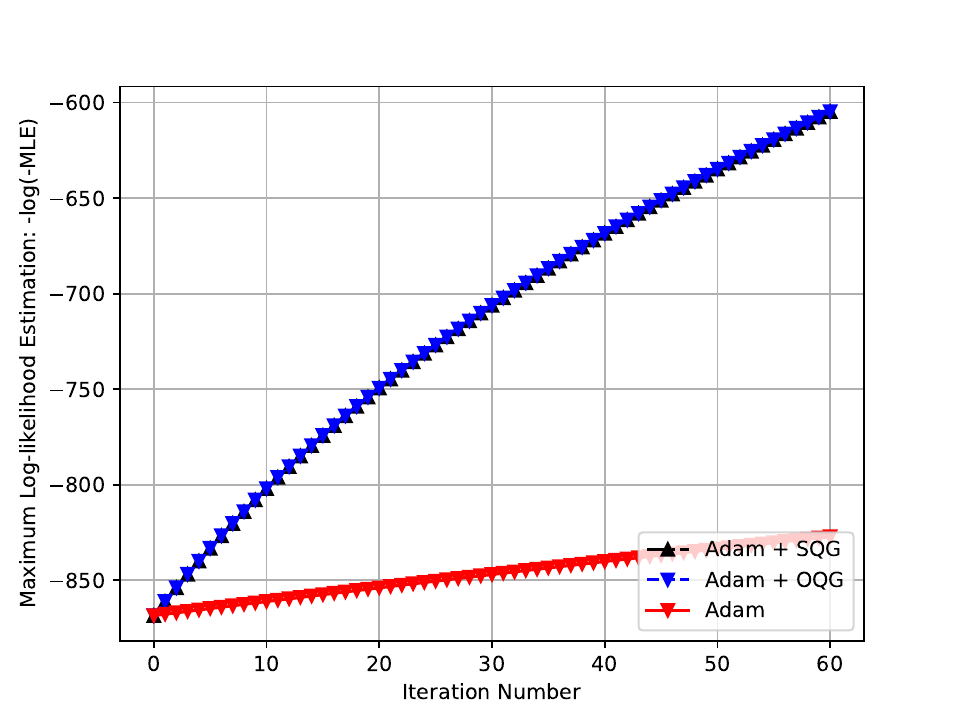}
    \label{fig:subfig02}
}

\vspace{1em} 

\subfloat[The lbw dataset]{%
    \includegraphics[width=0.48\textwidth]{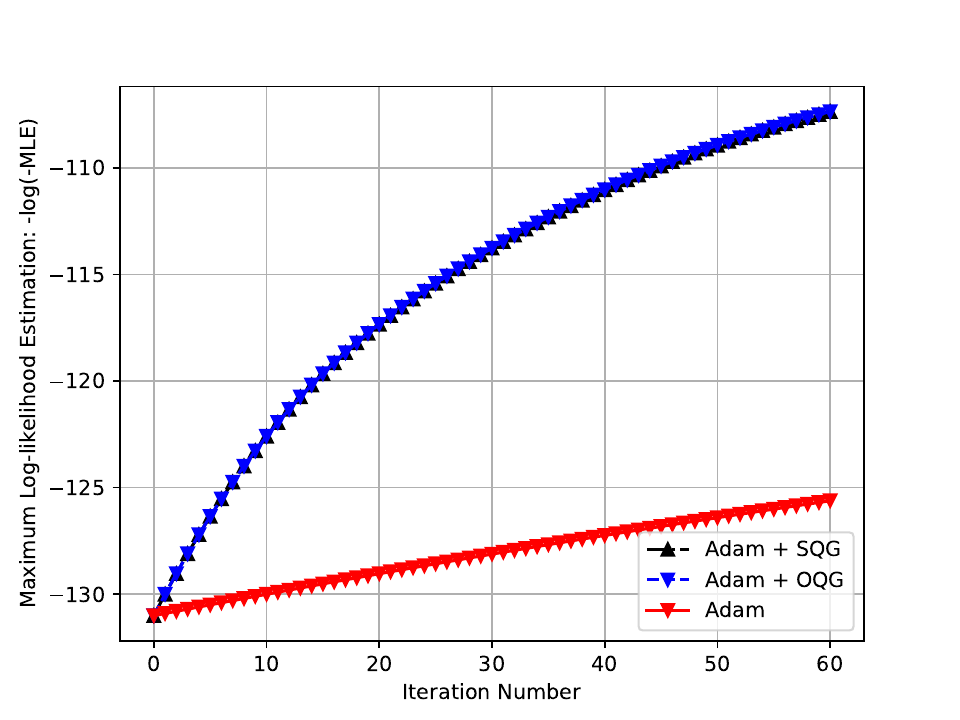}
    \label{fig:subfig03}
}
\hfill
\subfloat[The nhanes3 dataset]{%
    \includegraphics[width=0.48\textwidth]{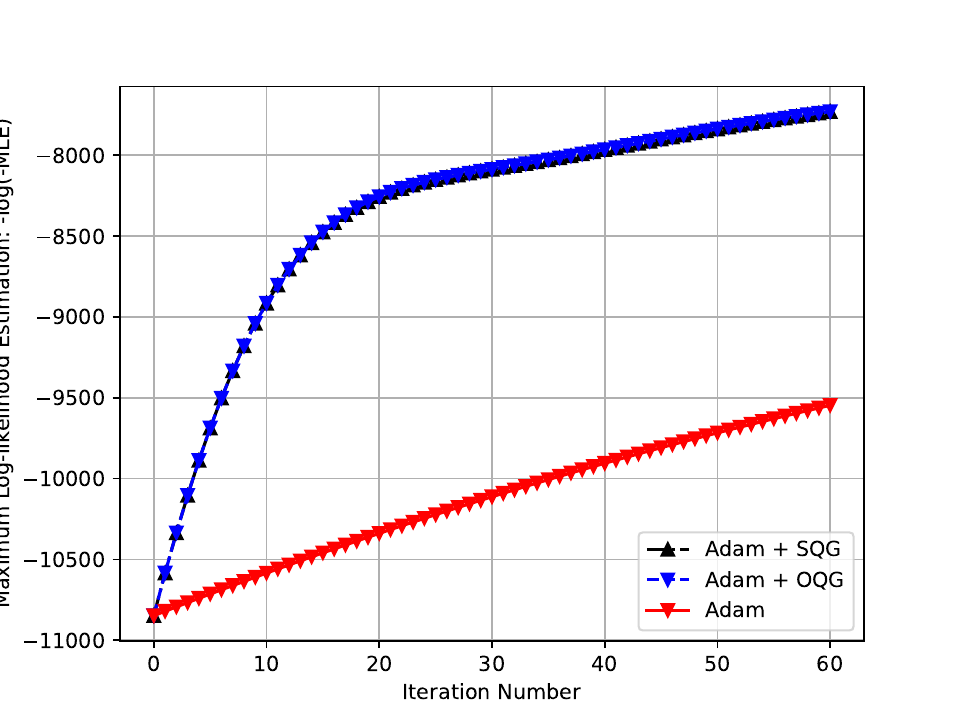}
    \label{fig:subfig04}
}

\vspace{1em} 

\subfloat[The pcs dataset]{%
    \includegraphics[width=0.48\textwidth]{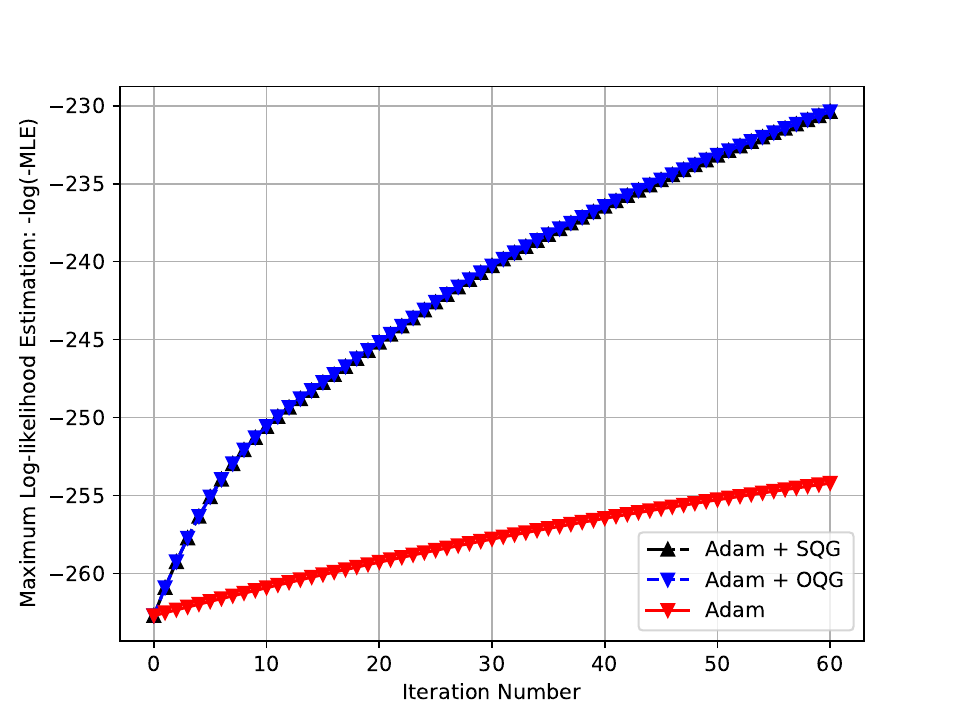}
    \label{fig:subfig03}
}
\hfill
\subfloat[The uis dataset]{%
    \includegraphics[width=0.48\textwidth]{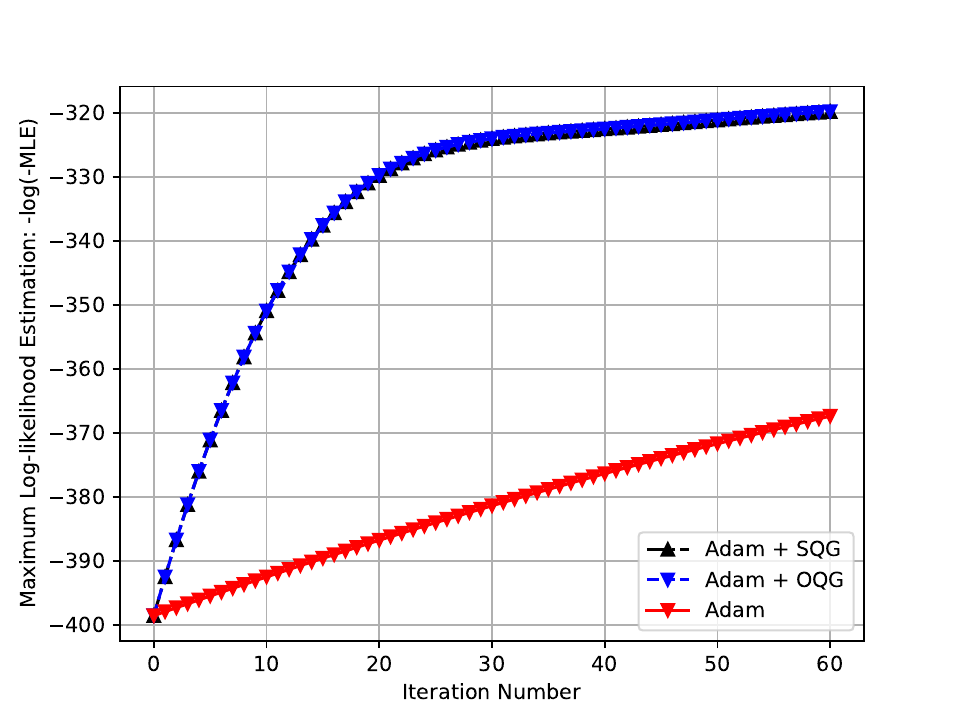}
    \label{fig:subfig04}
}

\vspace{1em} 

\subfloat[restructured MNIST dataset]{%
    \includegraphics[width=0.48\textwidth]{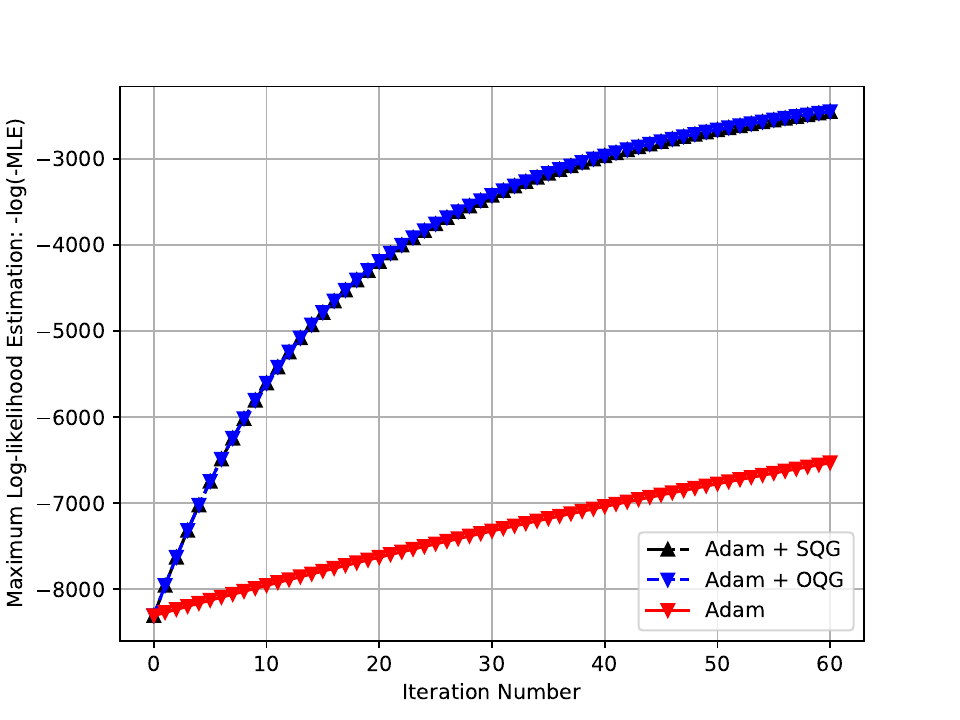}
    \label{fig:subfig03}
}
\hfill
\subfloat[The private financial dataset]{%
    \includegraphics[width=0.48\textwidth]{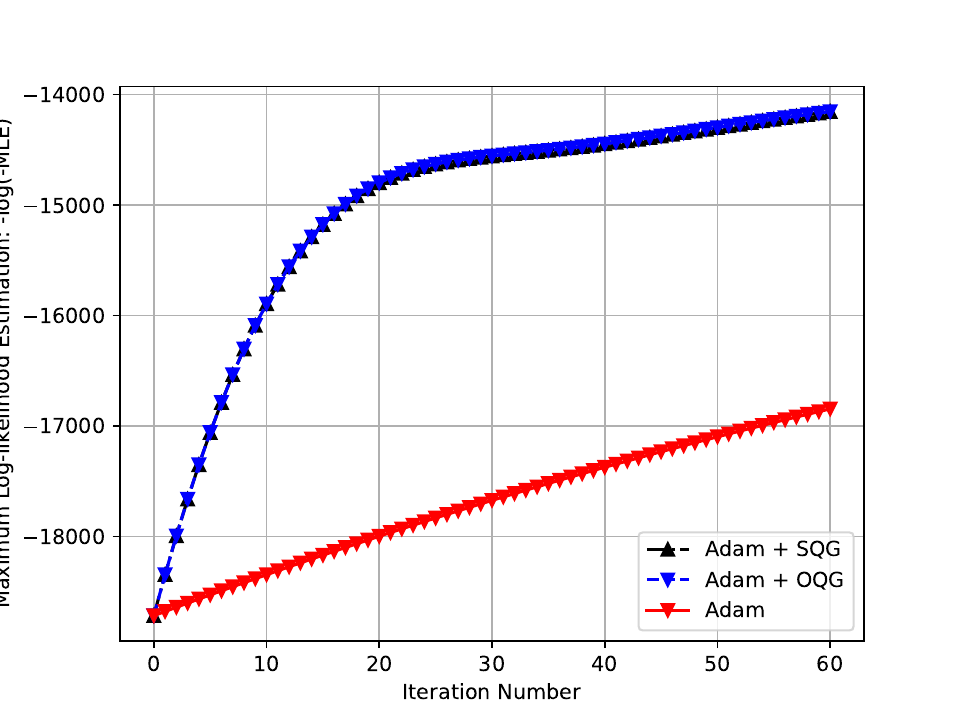}
    \label{fig:subfig04}
}

\caption{The training results of Adam + SQG vs. Adam + OQG vs. Adam in the clear domain.}
\label{fig2}
\end{figure}

\section{Conclusion}

In this work, we introduced \textit{Generalized Quadratic Gradient (GQG)}, a unified
framework that extends the quadratic gradient principle beyond specific Hessian
approximation strategies. By interpreting the quadratic gradient as the result of a
local quadratic model with a positive-definite curvature matrix, GQG provides a
general formulation for constructing curvature-aware optimization methods.

Unlike previous quadratic gradient approaches that rely on fixed Hessian approximations,
simplified diagonal approximations, or BFGS-based Hessian surrogates, GQG explores a
broader class of positive-definite curvature constructions. This perspective provides
a flexible foundation for extending quadratic-gradient-based optimization to various
Newton-type and quasi-Newton algorithms.

Future work will investigate more efficient positive-definite curvature matrix
constructions, explore theoretical convergence properties of generalized quadratic
gradient methods, and evaluate their applicability to large-scale machine learning
and distributed optimization.

All Python source code used to implement the experiments in this paper is openly
available at:
\href{https://github.com/petitioner/ML.GeneralizedQuadraticGradient}
{$\texttt{https://github.com/petitioner/ML.GeneralizedQuadraticGradient}$}.

Rather than designing a specific Hessian approximation, this work highlights that the
fundamental principle behind quadratic gradient optimization is the construction of an
appropriate positive-definite curvature geometry.

\bibliography{ML.GeneralizedQuadraticGradient}
\bibliographystyle{apalike}

\end{document}